\documentclass[12pt]{amsart}      %{article} was 12pt latex e
\usepackage{amssymb}
\usepackage{eucal}
\usepackage{amsmath}
\usepackage{amscd}
\usepackage[all]{xy}           %xypic macro for latex2.09

\usepackage{amsfonts,latexsym}
\usepackage{xspace}
\usepackage{hyperref}
\usepackage{float}
\usepackage{color}
\usepackage{colordvi}
\usepackage{multicol}

\newcommand{\nc}{\newcommand}
\newcommand{\delete}[1]{}

\nc{\dfootnote}[1]{{}}          %{{}}
\nc{\ffootnote}[1]{\dfootnote{#1}}

\delete{
\nc{\mfootnote}[1]{{}}        % Use this to suppress footnotes
\nc{\ofootnote}[1]{{}}        % Use this to suppress footnotes
}

\nc{\ofootnote}[1]{\footnote{\tiny Older version: #1}} % Use this to show footnotes

\nc{\mlabel}[1]{\label{#1}}  % Use this to suppress names
\nc{\mcite}[1]{\cite{#1}}  % Use this to suppress names
\nc{\mref}[1]{\ref{#1}}  % Use this to suppress names
\nc{\mkeep}[1]{{}}      % Use this to suppress marginpar
\nc{\mbibitem}[1]{\bibitem{#1}} % Use this to show number name
\delete{
\nc{\mcite}[1]{\cite{#1}{{\bf{{\ }(#1)}}}}  % Use this lines to show names
\nc{\mlabel}[1]{\label{#1}  % Use the next two lines to show names
{\hfill \hspace{1cm}{\bf{{\ }\hfill(#1)}}}}
\nc{\mref}[1]{\ref{#1}{{\bf{{\ }(#1)}}}}  % Use this lines to show names
\nc{\mbibitem}[1]{\bibitem[\bf #1]{#1}} % Use this to show name
\nc{\mkeep}[1]{\marginpar{{\bf #1}}} % Use this to show marginpar
}

\newtheorem{theorem}{Theorem}[section]
\newtheorem{prop}[theorem]{Proposition}
\newtheorem{defn}[theorem]{Definition}
\newtheorem{lemma}[theorem]{Lemma}
\newtheorem{coro}[theorem]{Corollary}
\newtheorem{prop-def}{Proposition-Definition}[section]

\nc{\bond}{\vdash}

\nc{\comp}[1]{\langle #1\rangle} \nc{\spr}{\cdot}
\nc{\disp}[1]{\displaystyle{#1}}
\nc{\sk}[1]{\mrm{sk}(#1)}
\nc{\ve}[1]{\mrm{vec}({#1})}
\nc{\Li}{\mrm{Li}}
\nc{\icut}{^!}
\nc{\bin}[2]{ (_{\stackrel{\scs{#1}}{\scs{#2}}})}  %binomial coeff
\nc{\binc}[2]{ \bigg (\!\! \begin{array}{c} \scs{#1}\\
    \scs{#2} \end{array}\!\! \bigg )}  %binomial coeff
\nc{\bbinc}[2]{ \left (\!\! \begin{array}{c} {#1}\\
    {#2} \end{array}\!\! \right )}  %binomial coeff
\nc{\bincc}[2]{  \left ( {\scs{#1} \atop
    \vspace{-.5cm}\scs{#2}} \right )}  %binomial coeff
\nc{\sarray}[2]{\begin{array}{c}#1 \vspace{.1cm}\\ \hline
    \vspace{-.35cm} \\ #2 \end{array}}
\nc{\spair}[2]{\big[\begin{array}{c}\scs{#1} \\ \scs{#2} \end{array} \big]}
\nc{\zsg}[1]{\widehat{#1}}
\nc{\ssg}[1]{\overline{#1}}
\nc{\bs}{\bar{S}}
\nc{\dcup}{\stackrel{\bullet}{\cup}}
\nc{\dbigcup}{\stackrel{\bullet}{\bigcup}} \nc{\la}{\longrightarrow}
\nc{\fe}{\'{e}} \nc{\rar}{\rightarrow} \nc{\dar}{\downarrow}
\nc{\dap}[1]{\downarrow \rlap{$\scriptstyle{#1}$}}
\nc{\uap}[1]{\uparrow \rlap{$\scriptstyle{#1}$}}
\nc{\dt}[1]{{#1}^\sharp} \nc{\st}[1]{{#1}^\flat}
\nc{\defeq}{\stackrel{\rm def}{=}} \nc{\dis}[1]{\displaystyle{#1}}
\nc{\dotcup}{\ \displaystyle{\bigcup^\bullet}\ } \nc{\hcm}{\
\hat{,}\ } \nc{\hcirc}{\hat{\circ}} \nc{\hts}{\hat{\shpr}}
\nc{\lts}{\stackrel{\leftarrow}{\shpr}}
\nc{\rts}{\stackrel{\rightarrow}{\shpr}} \nc{\lleft}{[}
\nc{\lright}{]} \nc{\uni}[1]{\tilde{#1}} \nc{\free}[1]{\bar{#1}}
\nc{\den}[1]{\check{#1}} \nc{\lrpa}{\wr} \nc{\curlyl}{\left \{
\begin{array}{c} {} \\ {} \end{array}
    \right . \!\!\!\!\!\!\!}
\nc{\curlyr}{ \!\!\!\!\!\!\!
    \left . \begin{array}{c} {} \\ {} \end{array}
    \right \} }
\nc{\longmid}{\left | \begin{array}{c} {} \\ {} \end{array}
    \right . \!\!\!\!\!\!\!}
\nc{\ot}{\otimes} \nc{\bigot}{\bigotimes} \nc{\mdiv}{\mrm{div}}
\nc{\sg}{G} \nc{\zg}{{Z}} \nc{\ig}{I} \nc{\pg}{P} \nc{\jg}{J}
\nc{\eg}{E} \nc{\fg}{F} \nc{\cg}{C} \nc{\mg}{M} \nc{\abg}{C}
\nc{\bas}{B} \nc{\Lyn}{\mrm{Lyn}} \nc{\lyn}{\Lyn} \nc{\sG}{{\cals}}
\nc{\iG}{{\cali}} \nc{\jG}{{\calj}} \nc{\eG}{{\cale}}
\nc{\pG}{{\calp}} \nc{\fG}{{\calf}} \nc{\cG}{{\calc}}
\nc{\mG}{{\calm}} \nc{\ora}[1]{\stackrel{#1}{\rar}}
\nc{\ola}[1]{\stackrel{#1}{\la}}%${\Bbb Z}$
\nc{\pex}[1]{\{#1\}} \nc{\scs}[1]{\scriptstyle{#1}}
\nc{\mrm}[1]{{\rm #1}} \nc{\sym}[1]{{\widehat{#1}}}
\nc{\margin}[1]{\marginpar{\rm #1}}   %{\rm #1}}
\nc{\dirlim}{\displaystyle{\lim_{\longrightarrow}}\,}
\nc{\invlim}{\displaystyle{\lim_{\longleftarrow}}\,}
\nc{\mvp}{\vspace{0.5cm}} \nc{\svp}{\vspace{2cm}}
\nc{\vp}{\vspace{8cm}} \nc{\proofbegin}{\noindent{\bf Proof: }}
\nc{\proofend}{$\blacksquare$ \vspace{0.5cm}}
\nc{\shqs}{\eta} \font\cyr=wncyr10
\nc{\sha}{{\mbox{\cyr X}}}  %used to be \cyr
\newfont{\scyr}{wncyr10 scaled 550}
\nc{\ssha}{\,\mbox{\bf \scyr X}\,}
\newfont{\bcyr}{wncyr10 scaled 1000}
\nc{\qssha}{{{\ssha\hspace{-2pt}_{\eta}}}\,}
\nc{\qsshae}{{{\ssha\hspace{-2pt}_{\eta}}}\,}
\nc{\qsshab}{{{\ssha\hspace{-2pt}_{\rho}}}\,}
\nc{\pssha}{{\star}}
\nc{\ncsha}{{\mbox{\cyr X}^{\mathrm NC}}} \nc{\ncshao}{{\mbox{\cyr
X}^{\mathrm NC,\,0}}}
\nc{\shpr}{\diamond}    %Shuffle product

\nc{\shf}{{^{\ssha}}} \nc{\qsh}{{^{\ast}}}
\nc{\bsh}{{^{\qsshab}}}
\nc{\esh}{{^{\qsshae}}}
\nc{\lshf}{_{\ssha}} \nc{\lqsh}{_{\ast}} \nc{\lzero}{_{\hskip -5pt
0}}
\nc{\shzero}{_{\hskip -7.5pt 0}}
\nc{\lizero}{{_{\hskip -10pt 0}}\,}
\nc{\lone}{_{\hskip -7.5pt 1}}
\nc{\shprl}{{{\shpr}_\lambda}}
\nc{\shpro}{\diamond^0}    %Shuffle product
\nc{\shpru}{\check{\diamond}} \nc{\catpr}{\diamond_l}
\nc{\rcatpr}{\diamond_r} \nc{\lapr}{\diamond_a}
\nc{\lepr}{\diamond_e} \nc{\tcon}{^{\ot}} \nc{\conv}{_c}
\nc{\vep}{\varepsilon} \nc{\labs}{\mid\!} \nc{\rabs}{\!\mid}
\nc{\hsha}{\widehat{\sha}} \nc{\lsha}{\stackrel{\leftarrow}{\sha}}
\nc{\rsha}{\stackrel{\rightarrow}{\sha}}
\nc{\EDS}{{\mrm{EDS}}\xspace} \nc{\DS}{{\mathbf{DS}}}
\nc{\lc}{[} \nc{\rc}{]} \nc{\rbset}{R} \nc{\rbnum}{r}
\nc{\rbfun}{\mathbf{R}} \nc{\pset}{P} \nc{\pnum}{p}
\nc{\pfun}{\mathbf{P}} \nc{\spset}{SP} \nc{\spnum}{sp}
\nc{\spgen}{\mathbf{SP}} \nc{\srbi}[1]{\{#1\}}

\nc{\rind}{r} \nc{\sind}{s} \nc{\tind}{t}  \nc{\kdim}{k}
\nc{\ldim}{\ell}

\nc{\ann}{\mrm{ann}} \nc{\Aut}{\mrm{Aut}} \nc{\can}{\mrm{can}}
\nc{\colim}{\mrm{colim}} \nc{\Cont}{\mrm{Cont}}
\nc{\rchar}{\mrm{char}} \nc{\cok}{\mrm{coker}} \nc{\dtf}{{R-{\rm
tf}}} \nc{\dtor}{{R-{\rm tor}}}

\nc{\Div}{{\mrm Div}} \nc{\End}{\mrm{End}} \nc{\Ext}{\mrm{Ext}}
\nc{\Fil}{\mrm{Fil}} \nc{\Frob}{\mrm{Frob}} \nc{\Gal}{\mrm{Gal}}
\nc{\GL}{\mrm{GL}} \nc{\lord}{\mrm{L-order}\xspace}
\nc{\rme}{\mrm{E}} \nc{\rmt}{\mrm{T}} \nc{\Sym}{\mrm{Sym}}
\nc{\Hom}{\mrm{Hom}} \nc{\hsr}{\mrm{H}} \nc{\hpol}{\mrm{HP}}
\nc{\id}{\mrm{id}} \nc{\im}{\mrm{im}} \nc{\incl}{\mrm{incl}}
\nc{\length}{\mrm{length}} \nc{\leng}{\mrm{\ell}} \nc{\LR}{\mrm{LR}}
\nc{\mchar}{\mrm char}
\nc{\MZV}{\mrm{MZV}\xspace}
\nc{\MZVs}{\mrm{MZVs}\xspace}
\nc{\MPV}{\mrm{MPV}\xspace}
\nc{\MPVs}{\mrm{MPVs}\xspace}
\nc{\MPL}{\mrm{MPL}\xspace}
\nc{\MPLs}{\mrm{MPLs}\xspace} \nc{\mzvalg}{\mathbf{MZV}}
\nc{\mplalg}{\mathbf{MPV}}
\nc{\edsalg}{\mathbf{EDS}} \nc{\qeds}{$\QQ$-EDS\xspace}
\nc{\zeds}{$\ZZ$-EDS\xspace} \nc{\zpeds}{$\ZZ_p$-EDS\xspace}
\nc{\fpeds}{$\FF_p$-EDS\xspace} \nc{\NC}{\mrm{NC}}
\nc{\mpart}{\mrm{part}} \nc{\os}{\mrm{OS}} \nc{\qs}{\mrm{QS}}
\nc{\ql}{{\QQ_\ell}} \nc{\qp}{{\QQ_p}} \nc{\rank}{\mrm{rank}}
\nc{\rcot}{\mrm{cot}} \nc{\rdef}{\mrm{def}} \nc{\rdiv}{{\rm div}}
\nc{\rtf}{{\rm tf}} \nc{\rtor}{{\rm tor}} \nc{\res}{\mrm{res}}
\nc{\sh}{\mrm{Sh}} \nc{\TL}{\mrm{TL}} \nc{\Spec}{\mrm{Spec}}
\nc{\tor}{\mrm{tor}} \nc{\Tr}{\mrm{Tr}} \nc{\tr}{\mrm{tr}}
\nc{\ETC}{\mathrm{ETC}} \nc{\ETL}{\mathrm{ETL}}
\nc{\EL}{\mathrm{EL}} \nc{\RETL}{\mathrm{RETL}}
\nc{\EETL}{\widetilde{\TL}} \nc{\word}{\rm word\xspace}
\nc{\words}{\rm words\xspace} \nc{\varab}{\phi_{\alpha,\beta}}

\nc{\ab}{\mathbf{Ab}} \nc{\Alg}{\mathbf{Alg}}
\nc{\Algo}{\mathbf{Alg}^0} \nc{\Bax}{\mathbf{Bax}}
\nc{\Baxo}{\mathbf{Bax}^0} \nc{\RBo}{\mathbf{RB}^0}
\nc{\BRB}{\mathbf{RB}} \nc{\Dend}{\mathbf{DD}} \nc{\bfe}{{\bf e}}
\nc{\bff}{{\bf f}} \nc{\bfk}{{\bf k}} \nc{\bfone}{{\bf 1}}
\nc{\base}[1]{{a_{#1}}} \nc{\detail}{\marginpar{\bf More detail}
    \noindent{\bf Need more detail!}
    \svp}
\nc{\Diff}{\mathbf{Diff}} \nc{\gap}{\marginpar{\bf
Incomplete}\noindent{\bf Incomplete!!}
    \svp}
\nc{\FMod}{\mathbf{FMod}} \nc{\RB}{\mathbf{RB}}
\nc{\Int}{\mathbf{Int}} \nc{\Mon}{\mathbf{Mon}}
\nc{\remarks}{\noindent{\bf Remarks: }} \nc{\Rep}{\mathbf{Rep}}
\nc{\Rings}{\mathbf{Rings}} \nc{\Sets}{\mathbf{Sets}}
\nc{\DT}{\mathbf{DT}}

\nc{\BA}{{\Bbb A}} \nc{\CC}{{\Bbb C}} \nc{\DD}{{\Bbb D}}
\nc{\EE}{{\Bbb E}} \nc{\FF}{{\Bbb F}} \nc{\GG}{{\Bbb G}}
\nc{\HH}{{\Bbb H}} \nc{\LL}{{\Bbb L}} \nc{\NN}{{\Bbb N}}
\nc{\QQ}{{\Bbb Q}} \nc{\RR}{{\Bbb R}} \nc{\TT}{{\Bbb T}}
\nc{\VV}{{\Bbb V}} \nc{\ZZ}{{\Bbb Z}}

\nc{\cala}{{\mathcal A}} \nc{\calc}{{\mathcal C}}
\nc{\cald}{{\mathcal D}} \nc{\cale}{{\mathcal E}}
\nc{\calf}{{\mathcal F}} \nc{\calg}{{\mathcal G}}
\nc{\calh}{{\mathcal H}} \nc{\cali}{{\mathcal I}}
\nc{\calj}{{\mathcal J}} \nc{\call}{{\mathcal L}}
\nc{\calm}{{\mathcal M}} \nc{\caln}{{\mathcal N}}
\nc{\calo}{{\mathcal O}} \nc{\calp}{{\mathcal P}}
\nc{\calr}{{\mathcal R}} \nc{\calt}{{\mathcal T}}
\nc{\calw}{{\mathcal W}} \nc{\calx}{{\mathcal X}}
\nc{\CA}{\mathcal{A}}

\nc\indI{\mathcal{I}}

\nc{\fraka}{{\mathfrak a}} \nc{\frakB}{{\mathfrak B}}
\nc{\frakb}{{\mathfrak b}} \nc{\frakd}{{\mathfrak d}}
\nc{\frakF}{{\mathfrak F}} \nc{\frakf}{{\mathfrak f}}
\nc{\frakg}{{\mathfrak g}} \nc{\frakL}{{\mathfrak L}}
\nc{\frakm}{{\mathfrak m}} \nc{\frakM}{{\mathfrak M}}
\nc{\frakMo}{{\mathfrak M}^0} \nc{\frakp}{{\mathfrak p}}
\nc{\fraks}{{\mathfrak s}}
\nc{\fraku}{{\mathfrak u}} \nc{\frakv}{{\mathfrak v}}
\nc{\frakw}{{\mathfrak w}}
\nc{\frakx}{{\mathfrak x}} \nc{\ox}{\overline{\frakx}}
\nc{\frakX}{{\mathfrak X}} \nc{\fraky}{{\mathfrak y}}
\nc{\frakz}{{\mathfrak z}}

\nc{\li}[1]{\textcolor{blue}{Li: #1}}
\nc{\wys}[1]{\textcolor{green}{William: #1}}
\nc{\byhs}[1]{\textcolor{red}{Bingyong: #1}}

\nc{\rrb}[1]{[#1]} \nc{\rrrb}[1]{\{#1\}} \nc{\ideal}[1]{\langle
#1\rangle} \nc{\refl}[1]{\overline{#1}} \nc{\rrB}{{reflexive
Rota-Baxter}\xspace} \nc{\rrrB}{{radical reflexive
Rota-Baxter}\xspace} \nc{\srb}{{strict Rota-Baxter}\xspace}
\nc{\Srb}{{Strict Rota-Baxter}\xspace}

\renewcommand\geq{\geqslant}
\renewcommand\leq{\leqslant}

\nc\rbop{{\lc\,\,\rc}} \nc\rbopi[1]{{\lc_{#1} \, \rc_{#1}}}

\nc{\redtext}[1]{\textcolor{red}{#1}}

\begin{document}

\title[Explicit double shuffle relations and Euler's formula]{Explicit double shuffle relations and a generalization of Euler's decomposition formula}
%
%========================================================================================%
\author{Li Guo}
\address{Department of Mathematics and Computer Science,
         Rutgers University,
         Newark, NJ 07102, USA}
\email{liguo@newark.rutgers.edu}
\author{Bingyong Xie}
\address{Department of Mathematics, Peking University, Beijing, 100871, China}
\email{byhsie@math.pku.edu.cn}

%=================================================================================
%\date{\today}
%==================================================================================
%\begin{document}

\begin{abstract}
We give an explicit formula for the shuffle relation in a general double shuffle framework that specializes to double shuffle relations of multiple zeta values and multiple polylogarithms. As an application, we generalize the well-known decomposition formula of Euler that expresses the product of two Riemann zeta values as a sum of double zeta values to a formula that expresses the product of two multiple polylogarithm values as a sum of other multiple polylogarithm values.
\medskip

\noindent
MSC classes: 11M41, 11M99, 40B05.
\smallskip

\noindent
Keywords: Euler's decomposition formula, multiple zeta values, multiple polylogarithm values, double shuffle relation.

\end{abstract}

\maketitle

%==================================================================================

%\tableofcontents

\setcounter{section}{0}

%================================================================================

\section {Introduction}
The decomposition formula of Euler is the equation
\begin{equation}
\zeta(r) \zeta(s)
 = \sum_{k=0}^{s-1} \binc{r+k-1}{k} \zeta(r+k,s-k)
+ \sum_{k=0}^{r-1} \binc{s+k-1}{k} \zeta(s+k,r-k), \quad r,s\geq 2,
\mlabel{eq:euler}
\end{equation}
expressing the product of two Riemann zeta values as a sum of double zeta values. In this paper we generalize this formula in two directions, from the product of one variable functions to that of multiple variables and from multiple zeta values to multiple polylogarithms. In fact, we obtain our formula in a general setting of shuffle algebras and quasi-shuffle algebras in order to provide a natural framework to treat these special values uniformly and to connect our generalization with the extended double shuffle relations of multiple zeta values.

To motivate our generalization, we explain the relationship between Euler's formula and the double shuffle relations of multiple zeta values.
Multiple zeta values (MZVs) have been studied quite intensively since the early 1990s~\mcite{Ho0,Za} involving many areas of mathematics and physics, from mixed Tate motives~\mcite{D-G,Te} to quantum field theory~\mcite{B-K}. Especially interesting are the algebraic and linear relations among the \MZVs. Because of the representations of an \MZV as an iterated sum and as an iterated integral, the multiplication of two \MZVs can be expressed in two ways as the sum of other \MZVs, one way following the {\bf quasi-shuffle (stuffle) relation} and the other way following the {\bf shuffle relation}. The combination of these two relations (called the {\bf double shuffle relations}) generates an extremely rich family of relations among \MZVs. In fact, as a conjecture, all relations among \MZVs can be derived from these relations and their degenerated forms, altogether called the {\bf extended double shuffle relations}~\mcite{IKZ,Ra}. A consequence of this conjecture is the irrationality of $\zeta(n)$ for all odd integers $n\geq 3$.

Naturally, determining  all the extended double shuffle relations is
challenging and the efforts have utilized a wide range
of methods. One difficulty is that the shuffle
relations have not been explicitly formulated in terms of the \MZVs.
For example, to determine the double shuffle relation from
multiplying two Riemann zeta values $\zeta(r)$ and $\zeta(s)$,
$r,s\geq 2$, one uses their sum representations and easily gets the
quasi-shuffle relation
\begin{equation}
 \zeta(r)\zeta(s)=\zeta(r,s)+\zeta(s,r)+\zeta(r+s).
\mlabel{eq:qsh2}
\end{equation}
On the other hand, to get their shuffle relation, one first  uses
their integral representations to express $\zeta(r)$ and $\zeta(s)$
as iterated integrals of dimensions $r$ and $s$, respectively. One
then uses the shuffle relation (or more concretely, repeated
applications of the integration by parts formula) to express the
product of these two iterated integrals as a sum of $\binc{r+s}{r}$
iterated integrals of dimension $r+s$. Finally, these last iterated
integrals are translated back to \MZVs and give the shuffle relation of $\zeta(r)\zeta(s)$. Explicitly, this shuffle relation is precisely the formula of Euler in Eq.~(\mref{eq:euler}). Then together with Eq.~(\mref{eq:qsh2}), we have the double
shuffle relation from $\zeta(r)$ and $\zeta(s)$.
See the recent papers~\mcite{BBG,3BL2,GKZ} for
the proofs of Euler's formula and see~\mcite{Br2,Zh} for its generalizations to double $q$-zeta values. For the applications of the double shuffle relation in this special case, we refer the reader to~\mcite{GKZ} on the connection of double zeta values with modular forms, and to~\mcite{OZ} on weighted sum formula of double zeta values.

In general, even though the computation of the shuffle relation can be performed recursively for any given pair of \MZVs, an explicit formula is missing so far.
As the above example shows, such an explicit formula not only provides an effective way to evaluate the shuffle relation, but also is important in the theoretical study of \MZVs, especially the double shuffle relations.
There are several families of
special values in addition to \MZVs, such as the alternating Euler sums~\mcite{3B}, the polylogarithms and
multiple polylogarithms~\mcite{3BL,Go}, especially at roots of unity~\mcite{Ra}, where the double shuffle relations are also studied~\mcite{BM,Ra,Zh3}, but are less understood.
Such an explicit formula for these values should also be useful to their study.

In this paper, we prove an explicit formula for the shuffle product in a general double shuffle framework. Consequently we obtain explicit shuffle formulas for the product of any two \MZVs, alternating Euler sums
and multiple polylogarithms, thereby generalizing Euler's formula.
As a concrete example, we obtain, for integers
$\rind_1,\sind_1\geq 2$ and $\sind_2\geq 1$,
\begin{equation}
\zeta(\rind_1)\,\zeta(\sind_1,\sind_2)  =
\hspace{-.6cm}\sum_{
\substack{\tind_1\geq 2,\tind_2\geq 1
\\ \tind_1+\tind_2=r_1+s_1}}   \hspace{-.4cm}\binc{\tind_1-1}{\rind_1-1}
\zeta(\tind_1,\tind_2,s_2)
\, +
\hspace{-.6cm}\sum_{
\substack{\tind_1\geq 2,\tind_2,\tind_3\geq 1
\\ \tind_1+\tind_2+\tind_3 \\
=\rind_1+\sind_1+\sind_2 }}
\hspace{-.2cm}
\binc{\tind_1-1}{\sind_1-1} \bigg
[\binc{\tind_2-1}{s_2-t_3}  +
\binc{\tind_2-1}{\sind_2-1} \bigg ]
\zeta(\tind_1,\tind_2,\tind_3).
\mlabel{eq:case1-2}
\end{equation}
As another instance, for integers $r_1,s_1\geq 2$ and $r_2,s_2\geq 1$, we have
{\allowdisplaybreaks
\begin{eqnarray}
\lefteqn{\zeta(\rind_1,\rind_2)\,\zeta(\sind_1,\sind_2)}
\notag
\\
=
 & \hspace{-.8cm}
\sum\limits_{ \substack{\tind_1\geq 2,\tind_2,\tind_3\geq 1
\\ \tind_1+\tind_2+\tind_3=r_1+r_2+s_1
}} \hspace{-.6cm}
 \binc{\tind_1-1}{\rind_1-1}\binc{\tind_2-1}{\rind_2-1}
\zeta(t_1,t_2,t_3,s_2)
 + \hspace{-.8cm} \sum\limits_{ \substack{\tind_1\geq 2,\tind_2,\tind_3\geq 1
\\ \tind_1+\tind_2+\tind_3=r_1+s_1+s_2
}} \hspace{-.6cm}
 \binc{\tind_1-1}{s_1-1}\binc{\tind_2-1}{s_2-1}
\zeta(t_1,t_2,t_3,r_2)
\notag \\
& \hspace{-3cm} +
\sum\limits_{ \substack{\tind_1\geq 2,\tind_2,\tind_3,\tind_4\geq 1
\\ \tind_1+\tind_2+\tind_3+\tind_4= \\
\rind_1+\rind_2+\sind_1+\sind_2
}}
\bigg [
\binc{\tind_1-1}{\rind_1-1} \binc{\tind_2-1}{t_1+t_2-r_1-s_1}
\bigg (\binc{\tind_3-1}{\sind_2-\tind_4}
 +  \binc{\tind_3-1}{\sind_2-1}\bigg)
\mlabel{eq:c2-2}
\\
& \quad \quad \quad\quad \quad  + \binc{\tind_1-1}{ \sind_1-1}\binc{\tind_2-1}{
t_1+t_2-r_1-s_1}\bigg (\binc{\tind_3-1}{\rind_2-\tind_4}
+\binc{\tind_3-1}{
\rind_2-1}\bigg) \bigg ]  \zeta(\tind_1,\tind_2,\tind_3,\tind_4).
\notag
\end{eqnarray}
}
We hope this framework can be further
extended to deal with other generalizations of multiple zeta values that have emerged recently,
such as the multiple $q$-zeta values~\mcite{Br,Zh} and renormalized \MZVs~\mcite{GZ,GZ2,MP}.

\smallskip

The organization of the paper is as follows.
In Section~\mref{sec:main}, we first describe the algebraic framework of double shuffle algebras. We then give our main formula in two variations (Theorem~\mref{thm:qshsh} and Theorem~\mref{thm:qshshe}). Theorem~\mref{thm:qshsh} is more general and easier to prove. Theorem~\mref{thm:qshshe} is more convenient for applications to multiple polylogarithm values and \MZVs (Corollary~\mref{co:mainmpl} and Corollary~\mref{co:mainmzv}). There we also provide some examples.
The proofs of the main theorems are quite long. So several lemmas are first proved in Section~\mref{sec:lem}. Then these lemmas are applied in Section~\mref{sec:proof} to prove the main formula by induction.
As an appendix, Section~\mref{sec:appen} includes a shuffle product formulation of the main formula.

\medskip

\noindent
{\bf Acknowledgements: } Both authors thank the hospitality and stimulating environment provided by the Max Planck Institute for Mathematics at Bonn where this research was carried out. They also thank Don Zagier and Matilde Marcolli for suggestions on an earlier draft and for encouragement.
The first author acknowledges the support from NSF grant DMS-0505643.

\section{The main theorems, applications and examples}
\mlabel{sec:main} We first set up in Section~\mref{ss:frame} a framework of general double shuffles to give a uniform formulation
of the double shuffle relations for multiple zeta values, alternating Euler sums and multiple polylogarithms. We then state in Section~\mref{ss:main} our main
formula in two variations in this framework.
Applications to the aforementioned special values are presented in Section~\mref{ss:app}. Computations in low dimensions and examples are provided in Section~\mref{ss:exam}.

\subsection{The general double shuffle framework}
\mlabel{ss:frame}
We formulate the framework to state our main theorems in Section~\mref{ss:main}. See Section~\mref{ss:app} for the concrete cases that have been considered before~\mcite{3BL,Go,Ho1,Ra,Zh3}.

We first introduce some notations. For any set $Y$, denote $M(Y)$ for the free monoid generated by $Y$.
Let $\calh(Y)$ be the free abelian group $\ZZ M(Y)$ with $M(Y)$ as a basis but without considering the product from the monoid $M(Y)$. When $\calh(Y)$ is equipped with an associative multiplication $\circ$, we use $\calh^{\circ}(Y)$ to denote the algebra $(\calh(Y),\circ)$.

Let $\sg$ be a given set. Define
$$\ssg{\sg}=\{x_0\}\cup\{x_b\ |\ b\in \sg\}$$
to be the set of symbols indexed by the disjoint set $\{0\}\sqcup \sg.$
Then the shuffle algebra~\mcite{Ka,Re} generated by $\ssg{\sg}$ is
\begin{equation}
\calh\shf(\ssg{\sg}):= (\ZZ M (\ssg{\sg}), \ssha)
\mlabel{eq:shalg}
\end{equation}
where the shuffle product $\ssha$ is defined recursively by
$$ (a_1 \fraka) \ssha (b_1 \frakb) = a_1  (\fraka \ssha (b_1  \frakb))+b_1  ((a_1 \fraka)\ssha \frakb), a_1,b_1\in \ssg{\sg}, \fraka,\frakb \in M (\ssg{\sg})$$ with the convention that $1\ssha \frakb=\frakb=\frakb \ssha 1$ for $\frakb \in
M(\ssg{\sg})$.
Define the subalgebra
\begin{equation}
\calh\shf\lone(\ssg{\sg}):=\ZZ \oplus \big(\oplus_{b\in \sg}
\calh\shf(\ssg{\sg}) x_b\big) \subseteq \calh\shf(\ssg{\sg}).
\mlabel{eq:shsubg}
\end{equation}

For the given set $\sg$, let $\zsg{\sg}$ be the set product
$$
\zsg{\sg}:=\ZZ_{\geq 1}\times \sg=\{w:=\spair{s}{b}\ |\ s \in
\ZZ_{\geq 1},
    b\in \sg\}.
$$
We will denote the non-unit elements in the free monoid $M(\zsg{\sg})$ by vectors
$$ \vec{\nu}:=[\nu_1,\cdots,\nu_k]
 =\spair{\sind_1,\cdots,\sind_k}{b_1,\cdots,b_k} =\spair{\vec{\sind}}{\vec{b}}\in \zsg{\sg}^k$$
and denote $[\nu_1,[\nu_2,\cdots,\nu_k]]=[\nu_1,\nu_2,\cdots,\nu_k].$
Consider the free abelian group
$$\calh(\zsg{\sg}):= \ZZ M(\zsg{\sg})=\bigoplus_{\vec{\nu}\in \zsg{\sg}^k,\, k\geq 0}
\ZZ \vec{\nu}, \quad \zsg{\sg}^0=\{1\}. $$

As in the case of the shuffle algebra from \MZVs, elements of $\calh\shf\lone(\ssg{\sg})$ of the form $$x_0^{s_1-1}x_{b_1}x_0^{s_2-1}x_{b_2}\cdots x_0^{s_k-1}x_{b_k}, \quad s_i\geq 1, b_i\in G, 1\leq i\leq k, k\geq 1,$$
together with $1$, form a basis of $\calh\shf\lone(\ssg{\sg}).$ Since $\calh(\zsg{\sg})$ with the concatenation product is the free non-commutative algebra generated by $\zsg{\sg}$, there is a natural linear bijection
\begin{equation}
\rho: \calh\shf\lone(\ssg{\sg}) \to \calh(\zsg{\sg}), \quad
 x_0^{\sind_1-1}  x_{b_1}  \cdots   x_0^{\sind_k-1}  x_{b_k}
 \leftrightarrow \spair{ \sind_1,\  \sind_2,\ \cdots, \ \sind_k}{b_1,\ b_2,\ \cdots,\ b_k}, \quad 1\leftrightarrow 1.
\mlabel{eq:beta}
\end{equation}
Through $\rho$, the shuffle product $\ssha$ on $\calh\shf\lone(\ssg{\sg})$ defines a product on $\calh(\zsg{\sg})$ by
\begin{equation}
\vec{\mu}\qsshab \vec{\nu}: = \rho( \rho^{-1}(\vec{\mu})\ssha
\rho^{-1}(\vec{\nu})), \quad \vec{\mu},\vec{\nu}\in
\calh(\zsg{\sg}). \mlabel{eq:shtransb}
\end{equation}
Following our notations, we use $\calh\bsh(\zsg{\sg})$ to denote this algebra.

Now assume that $\sg$ is a multiplicative abelian group.
Equip $\zsg{\sg}=\ZZ_{\geq 1}\times \sg$ with the abelian semigroup structure by the component multiplication: $\spair{\sind_1}{z_1}\cdot
\spair{\sind_2}{z_2}=\spair{\sind_1+\sind_2}{z_1z_2}.$
Then we define the quasi-shuffle algebra~\mcite{Ho2} on $\zsg{\sg}$ to be
\begin{equation}
\calh\qsh(\zsg{\sg}):= (\ZZ M(\zsg{\sg}),\ast)
\mlabel{eq:qshalg}
\end{equation}
where the multiplication $\ast$ is defined by the recursion
$$ [\mu_1,\vec{\mu}\,'] \ast [\nu_1,\vec{\nu}\,'] = [\mu_1,(\vec{\mu}\,'\ast [\nu_1,\vec{\nu}\,'])]
+ [\nu_1,[\mu_1, \vec{\mu}\,']\ast \vec{\nu}\,'] + [(\mu_1\cdot
\nu_1),\vec{\mu}\,'\ast \vec{\nu}\,'], $$
$\mu_1,\nu_1\in \zsg{\sg},
\vec{\mu}\,',\vec{\nu}\,'\in M(\zsg{\sg}),
$
with the initial condition that $1\ast \vec{\nu}=\vec{\nu}=\vec{\nu}\ast 1$
for $\vec{\nu}\in M(\zsg{\sg})$. See~\mcite{G-K1,GX,Ho2} for its explicit description and its structure.

We define a linear bijection
\begin{equation}
\theta: \calh\qsh(\zsg{\sg})  \to \calh\qsh(\zsg{\sg}),
\quad
    \spair{s_1,\cdots,s_k}{b_1,\cdots,b_k}
     \mapsto \spair{s_1,\,s_2,\cdots,\,s_k}{\frac{1}{b_1},\,\frac{b_1}{b_2}, \cdots,\, \frac{b_{k-1}}{b_k}}
\mlabel{eq:la}
\end{equation}
whose inverse is given by
\begin{equation}
\theta^{-1}: \calh\qsh(\zsg{\sg})  \to \calh\qsh(\zsg{\sg}), \quad
    \spair{s_1,\cdots,s_k}{z_1,\cdots,z_k}
     \mapsto \spair{s_1,\,s_2,\cdots,\,s_k}{\frac{1}{z_1},\,\frac{1}{z_1z_2}, \cdots,\, \frac{1}{z_1\cdots z_k}}
\mlabel{eq:lainv}
\end{equation}
Note that the action of $\theta$ is defined by an action on the second row of elements in $\calh\qsh(\zsg{\sg})$ which is again denoted by $\theta$:
\begin{equation}
\theta(b_1,\cdots,b_k) = \big(\frac{1}{b_1},\,\frac{b_1}{b_2}, \cdots,\, \frac{b_{k-1}}{b_k}\big).
\mlabel{eq:la1}
\end{equation}

The composition of $\rho$ and $\theta$ gives a natural bijection of abelian groups (but {\em not} as algebras)
\begin{equation}
\shqs: \calh\shf\lone(\ssg{\sg}) \to \calh\qsh(\zsg{\sg}), \quad
 x_0^{\sind_1-1}  x_{b_1}  \cdots   x_0^{\sind_k-1}  x_{b_k}
 \leftrightarrow \spair{ \sind_1,\  \sind_2,\ \cdots, \ \sind_k}{\frac{1}{b_1},\ \frac{b_1}{b_2},\ \cdots,\ \frac{b_{k-1}}{b_k}}
\mlabel{eq:shqsh1}
\end{equation}
whose inverse is given by
$\spair{\sind_1,\cdots,s_k}{z_1,\cdots,z_k} \mapsto
x_0^{s_1-1}x_{z_1^{-1}}x_0^{s_2-1}x_{(z_1z_2)^{-1}} \cdots
x_0^{s_k-1}x_{(z_1\cdots z_k)^{-1}}.$

Through $\shqs$, the shuffle product $\ssha$ on
$\calh\shf\lone(\ssg{\sg})$ transports to a product $\qssha$ on
$\calh(\zsg{\sg})$, resulting a commutative algebra $\calh\esh(\zsg{\sg})=(\calh(\zsg{\sg}), \qssha)$. More precisely, for $\vec{\mu},\vec{\nu}\in
\calh(\zsg{\sg})$,
\begin{equation}
\vec{\mu}\qssha \vec{\nu}: = \shqs( \shqs^{-1}(\vec{\mu})\ssha
\shqs^{-1}(\vec{\nu})). \mlabel{eq:shtrans}
\end{equation}
Then we have the following commutative diagram of commutative algebras:
\begin{equation}
\xymatrix{
\calh\shf\shzero(\ssg{S^1}) \ar^{\rho}[rr]  \ar@/^2.5pc/^{\eta}[rrrr]&&
    \calh\bsh\lizero(\zsg{S}^1)
    \ar^{\theta}[rr] && \calh\esh\lizero(\zsg{S}^1)
}
\mlabel{eq:comm}
\end{equation}
The purpose of this paper is to give an explicit formula for
$\vec{\mu} \qsshae \vec{\nu}$ (Theorem~\mref{thm:qshshe}) which naturally gives shuffle formulas
for \MZVs, \MPVs and alternating Euler sums.
However, as we will see later, for the proof of this formula, it is more convenient to work with its variation (Theorem~\mref{thm:qshsh}) for the product $\qsshab$ since it is more compatible with the module structure on $\calh\qsh(\zsg{\sg})$. This approach also allows us to obtain a formula without requiring that $G$ is a group, further extending its potential of applications that will be discussed in a future work.

\subsection{The statement of the main theorems} \mlabel{ss:main}
We first introduce some notations.
For positive integers $k$ and $\ell$, denote
$[k]=\{1,\cdots,k\}$ and $[k+1,k+\ell]=\{k+1,\cdots,k+\ell\}.$
Define
\begin{equation}
\indI_{k,\ell}=\left \{(\varphi,\psi)\ \Big|\ \begin{array}{l}
\varphi: [k]\to [k+\ell], \psi: [\ell]\to [k+\ell]
\text{ are order preserving } \\
\text{ injective maps and } \im(\varphi)\sqcup\im (\psi)=[k+\ell] \end{array} \right
\} \mlabel{eq:ind}
\end{equation}
Let $\vec{a}\in\sg^k$, $\vec{b}\in\sg^\ell$ and
$(\varphi,\psi)\in\indI_{k,\ell}$. We define
$\vec{a}\ssha_{(\varphi,\psi)}\vec{b}$ to be the vector whose $i$th
component is
\begin{equation}
(\vec{a}\ssha_{(\varphi,\psi)} \vec{b})_i =\left\{\begin{array}{ll}
a_j & \text{if } i=\varphi(j) \\ b_j & \text{if }
i=\psi(j)\end{array}\right.
= a_{\varphi^{-1}(i)}b_{\psi^{-1}(i)}, \quad 1\leq i\leq k+\ell,
\mlabel{eq:mulind}
\end{equation}
with the convention that $a_\emptyset=b_\emptyset=1.$

Let
$\vec{\rind}=(\rind_1,\cdots, \rind_k)\in\ZZ_{\geq 1}^k$, $\vec{\sind}=(\sind_1,\cdots,\sind_\ell)\in \ZZ_{\geq 1}^{\ell}$ and  $\vec{\tind}=(\tind_1,\cdots, \tind_{k+\ell})\in \ZZ_{\geq 1}^{k+\ell}$
with $|\vec{\rind}|+|\vec{\sind}|=|\vec{\tind}|$. Here
$|\vec{\rind}|=\rind_1+\cdots +\rind_k$ and similarly for $|\vec{\sind}|$ and $|\vec{\tind}|$.
Denote $R_i=r_1+\cdots +r_i$ for $i\in [k]$, $S_i=s_1+\cdots +s_i$ for $i\in [\ell]$ and $T_i=t_1+\cdots+t_i$ for $i\in [k+\ell]$.
For $i\in [k+\ell]$, define
\begin{equation}
h_{(\varphi,\psi),i}=h_{(\varphi,\psi),(\vec{\rind},\vec{\sind}),i}
=
       \left\{
              \begin{array}{ll} \rind_{j} & \text{ if } i=\varphi(j)
                                 \\
                                \sind_{j} & \text{ if } i=\psi(j)
              \end{array}
              \right.
              = r_{\varphi^{-1}(i)}s_{\psi^{-1}(i)},
\mlabel{eq:h0}
\end{equation}
with the convention that $r_\emptyset =s_\emptyset =1.$

With these notations, we define
\begin{equation}
c_{\vec{\rind},\vec{\sind}}^{\vec{\tind},(\varphi,\psi)}(i) =\left\{
\begin{array}{ll}
\binc{\tind_i-1}{h_{(\varphi,\psi),i}-1} &
    \begin{array}{l}\text{if } i=1,  \text{if }
i-1,i \in \im(\varphi)\\ \text{or if }
i-1,i \in \im(\psi),
\end{array}
\vspace{.2cm}
\\ \vspace{.2cm}
\begin{array}{l}
\binc{\tind_i-1} {T_i-R_{|\varphi^{-1}([i])|}-S_{|\psi^{-1}([i])|}}\\
= \binc{\tind_{i}-1}{\sum\limits_{j=1}^{i} \tind_j
-\sum\limits_{j=1}^{i} h_{(\varphi,\psi),j}}
\end{array} & \text{ otherwise}.
\end{array}
\right. \mlabel{eq:coef-re-def1}
\end{equation}
Denote \begin{equation} c_{\vec{\rind},\vec{\sind}}^{\vec{\tind},
(\varphi,\psi)}
:=\prod_{i=1}^{k+\ell}c_{\vec{\rind},\vec{\sind}}^{\vec{\tind},
(\varphi,\psi)}(i)=\prod_{j=1}^k
c_{\vec{\rind},\vec{\sind}}^{\vec{\tind},
(\varphi,\psi)}(\varphi(j))\prod_{j=1}^\ell
c_{\vec{\rind},\vec{\sind}}^{\vec{\tind}, (\varphi,\psi)}(\psi(j)).
\mlabel{eq:coef-sh-eq}
\end{equation}

Now we can state the first variation of our main formula.

\begin{theorem} %{\bf (Main Theorem: the Algebraic Form) }
Let $k,\ell$ be positive integers. Let $\sg$ be a set and let $\calh\bsh(\zsg{\sg})=(\calh(\zsg{\sg}), \qsshab)$ be the algebra defined by Eq.~(\mref{eq:shtransb}).
Then for $\spair{\vec{\rind}}{\vec{a}}\in \zsg{\sg}^k$ and
$\spair{\vec{\sind}}{\vec{b}}\in\zsg{\sg}^\ell$ in $\calh\bsh(\zsg{\sg})$, we have
\begin{equation}
\begin{aligned}
  \spair{\vec{\rind}}{\vec{a}}\qsshab \spair{\vec{\sind}}{\vec{b}}
  &=
\sum_{\substack{ (\varphi,\psi)\in \indI_{k,\ell}\\
\vec{\tind}\in \ZZ_{\geq 1}^{k+\ell},
|\vec{\tind}|=|\vec{\rind}|+|\vec{\sind}| }}
c_{\vec{\rind},\vec{\sind}}^{\vec{\tind},
  (\varphi,\psi)}
  \spair{\vec{\tind}}{\vec{a}\ssha_{(\varphi,\psi)}\vec{b}}
  \\
  &=
\sum_{\substack{ (\varphi,\psi)\in \indI_{k,\ell}\\
\vec{\tind}\in \ZZ_{\geq 1}^{k+\ell},
|\vec{\tind}|=|\vec{\rind}|+|\vec{\sind}| }}
\bigg(\prod_{i=1}^{k+\ell}c_{\vec{\rind},\vec{\sind}}^{\vec{\tind},
(\varphi,\psi)}(i)\bigg)
\spair{\vec{\tind}}{\vec{a}\ssha_{(\varphi,\psi)}\vec{b}},
\end{aligned} \mlabel{eqn:maincoef}
\end{equation}
where $c_{\vec{\rind},\vec{\sind}}^{\vec{\tind}, (\varphi,\psi)}(i)$
is given in Eq.~(\mref{eq:coef-re-def1}) and
$\vec{a}\ssha_{(\varphi,\psi)}\vec{b}$ is given in
Eq.~(\mref{eq:mulind}). \mlabel{thm:qshsh}
\end{theorem}

For the purpose of applications to \MZVs and multiple polylogarithms, we give an equivalent form of Theorem \mref{thm:qshsh} under the condition that $\sg$ is an abelian group. For
$\vec{w}\in \sg^k$ and $\vec{z}\in \sg^\ell$, we define
\begin{equation}
(\vec{w}\pssha_{(\varphi,\psi)}\vec{z})_i =\left\{\begin{array}{ll}
w_j & \text{ if } i=\varphi(j) \text{ and } \text{ either } i=1 \text{ or } i-1\in\im(\varphi),   \smallskip \\
z_j & \text{ if } i=\psi(j) \text{ and } \text{ either } i=1 \text{ or } i-1\in\im(\varphi),  \smallskip \\
\frac{w_1\cdots w_j}{z_1\cdots z_{i-j}} & \text{ if } i=\varphi(j)
\text{ and } i-1\in \im(\psi), \smallskip\\
\frac{z_1\cdots z_j}{w_1\cdots w_{i-j}} & \text{ if } i=\psi(j)
\text{ and } i-1\in \im(\varphi).
\end{array}\right. \mlabel{eq:pssha}
\end{equation}

\begin{theorem} %{\bf (Main Theorem: the Arithmetic Form) }
Let $k,\ell$ be positive integers. Let $\sg$ be an abelian group and let $\calh\esh(\zsg{\sg})=(\calh(\zsg{\sg}), \qsshae)$ be the algebra defined by  Eq.~(\mref{eq:shtrans}).
Then for $\spair{\vec{\rind}}{\vec{w}}\in \zsg{\sg}^k$ and
$\spair{\vec{\sind}}{\vec{z}}\in\zsg{\sg}^\ell$ in $\calh\esh(\zsg{\sg})$, we have
\begin{equation}
\begin{aligned}
  \spair{\vec{\rind}}{\vec{w}}\qsshae \spair{\vec{\sind}}{\vec{z}}
  &=
\sum_{\substack{ (\varphi,\psi)\in \indI_{k,\ell}\\
\vec{\tind}\in \ZZ_{\geq 1}^{k+\ell},
|\vec{\tind}|=|\vec{\rind}|+|\vec{\sind}| }}
c_{\vec{\rind},\vec{\sind}}^{\vec{\tind},
  (\varphi,\psi)}
  \spair{\vec{\tind}}{\vec{w}\pssha_{(\varphi,\psi)}\vec{z}}
  \\
  &=
\sum_{\substack{ (\varphi,\psi)\in \indI_{k,\ell}\\
\vec{\tind}\in \ZZ_{\geq 1}^{k+\ell},
|\vec{\tind}|=|\vec{\rind}|+|\vec{\sind}| }}
\bigg(\prod_{i=1}^{k+\ell}c_{\vec{\rind},\vec{\sind}}^{\vec{\tind},
(\varphi,\psi)}(i)\bigg)
\spair{\vec{\tind}}{\vec{w}\pssha_{(\varphi,\psi)}\vec{z}},
\end{aligned} \mlabel{eqn:maincoefe}
\end{equation}
where $c_{\vec{\rind},\vec{\sind}}^{\vec{\tind}, (\varphi,\psi)}(i)$
is given in Eq.~(\mref{eq:coef-re-def1}) and
$\vec{w}\pssha_{(\varphi,\psi)}\vec{z}$ is given in
Eq.~(\mref{eq:pssha}). \mlabel{thm:qshshe}
\end{theorem}

We will next give applications and examples of Theorem
\mref{thm:qshshe} in Section~\mref{ss:app} and
Section~\mref{ss:exam}. Theorem~\mref{thm:qshshe} will be shown to follow from Theorem~\mref{thm:qshsh} in Section~\mref{ss:equiv}, and Theorem~\mref{thm:qshsh} will be proved in Section~\mref{ss:proof}. Preparational lemmas will be given in Section~\mref{sec:lem}.

\subsection{Applications}
\mlabel{ss:app} In this section, Theorem~\mref{thm:qshshe} is
specialized to give formulas for multiple zeta values, alternating
Euler sums and multiple polylogarithms. We start with multiple
polylogarithms and then specialize further to \MZVs and alternating
Euler sums. In Section~\mref{ss:exam} we demonstrate how to apply
these formulas by computing examples in low dimensions.

\subsubsection{Multiple polylogarithms}
\mlabel{sss:mpl}
A {\bf multiple polylogarithm value (\MPV)}~\mcite{3BL,Go,Go2} is defined by
\begin{equation}
\Li_{s_1,\cdots,s_k}(z_1,\cdots,z_k):= \sum_{n_1>\cdots>n_k\geq 1}
\frac{z_1^{n_1}\cdots z_k^{n_k}}{n_1^{s_1}\cdots n_k^{s_k}}
\mlabel{eq:mpls}
\end{equation}
where $|z_i|\leq 1$, $s_i\in \ZZ_{\geq 1}$, $1\leq i\leq k$, and $(s_1,z_1)\neq (1,1)$.
When $z_i=1, 1\leq i\leq k$, we obtain the multiple zeta values
$
 \zeta(s_1,\cdots,s_k)$
that we will consider in Section~\mref{sss:mzv}.
More generally, the special cases when $z_i$ are roots of unity have been studied~\mcite{3BL,BB,Go2,Ra} in connection with high cyclotomic theory, mixed motives and combinatorics, and have been found in the computations of Feynman diagrams~\mcite{Bro}.

With the notation of~\mcite{3BL}, we have
\begin{equation}
\begin{aligned}
&\Li_{s_1,\cdots,s_k}(z_1,\cdots,z_k)= \lambda\big(\begin{array}{c} s_1,\cdots,s_k\\ b_1,\cdots,b_k \end{array} \big):=
\sum_{n_1>n\cdots >n_k\geq 1}\frac{\big(\frac{1}{b_1}\big)^{n_1} \big(\frac{b_1}{b_2}\big)^{n_2}\cdots \big(\frac{b_{k-1}}{b_k}\big)^{n_k}} {n_1^{s_1}n_2^{s_2}\cdots n_k^{s_k}}, \\
&  \text{where } (b_1,\cdots,b_k)= \theta^{-1}(z_1,\cdots,z_k)=(z_1^{-1}, (z_1z_2)^{-1}, \cdots,
(z_1\cdots z_k)^{-1}).
\end{aligned}
\mlabel{eq:3b}
\end{equation}
Here $\theta$ is as defined in Eq.~(\mref{eq:la1}).

The product of two sums representing two \MPVs is a $\ZZ$-linear combination of other such sums. So the $\ZZ$-linear span of these values is an algebra which we denote by
$$
\mplalg= \ZZ \{ \Li_{s_1,\cdots,s_k}(z_1,\cdots,z_k)\ |\ s_i\in \ZZ_{\geq 1}, |z_i|\leq 1, (s_1,z_1)\neq (1,1)\}.
$$

In the framework of Section~\mref{ss:frame} and Section~\mref{ss:main}, take $\sg$ to be the multiplicative abelian group $S^1:=\{z\in \CC^{\times}\ |\ |z|= 1\}$, and consider the subalgebra
$$
\calh\qsh\lzero(\zsg{S}^1):=\ZZ \oplus \Big (\bigoplus_{\tiny\spair{s_1}{z_1}\neq \spair{1}{1}} \ZZ
\spair{s_1,s_2,\cdots, s_k}{z_1,z_2,\cdots,z_k} \Big )
\subseteq \calh\qsh(\zsg{S}^1). $$
Then $\calh\qsh(\zsg{\sg})$ coincides with the quasi-shuffle (stuffle) algebra~\mcite{Go2,Ra}  encoding \MPVs, and the multiplication rule of two \MPVs according to their sum representations in Eq.~(\mref{eq:mpls}) follows from the fact that the linear map
$$
\Li\qsh: \calh\qsh\lzero(\zsg{S}^1) \to \mplalg,
\quad \spair{s_1,\cdots,s_k}{z_1,\cdots,z_k} \mapsto \Li_{s_1,\cdots,s_k}(z_1,\cdots,z_k)
$$
is an algebra homomorphism.

We also consider the shuffle algebra $\calh\shf(\ssg{S^1 })$ and its subalgebras
$$
\begin{aligned}
 \calh\shf\shzero(\ssg{S^1}) &:= \ZZ \oplus \big(\oplus_{a, b\in \{0\}\cup S^1, a\neq 1, b\neq 0 } x_a\calh\shf(\ssg{S^1}) x_b\big)\\
& \subseteq \calh\shf\lone(\ssg{S^1}):=\ZZ \oplus \big(\oplus_{b\in
S^1 } \calh\shf(\ssg{S^1}) x_b\big) \subseteq \calh\shf(\ssg{S^1}).
\end{aligned}
$$
They agree with the shuffle algebras~\mcite{Go,Ra} encoding \MPVs through their integral representations~\mcite{3BL,Go,Ra}
\begin{equation}
\Li_{s_1,\cdots,s_k}(z_1,\cdots,z_k) =
 \int_0^1 \int_0^{u_1}\cdots \int_0^{u_{|\vec{s}|-1}} \frac{du_1}{f_1(u_1)}\cdots \frac{du_{|\vec{s}|}}{f_{|\vec{s}|}(u_{|\vec{s}|})}.
 \mlabel{eq:intmpl}
\end{equation}
Here
$$f_j(u_j)=\left\{\begin{array}{ll} (z_1\cdots z_i)^{-1}-u_j& \text{if } j= s_1+\cdots +s_i, 1\leq i \leq k,\\
u_j & \text{otherwise}. \end{array} \right .
$$

It takes a simpler form in terms of $\lambda\big(\begin{array}{c} s_1,\cdots,s_k\\ b_1,\cdots,b_k \end{array} \big)$ thanks to Eq.~(\mref{eq:3b}):
\begin{equation}
\lambda\big(\begin{array}{c} s_1,\cdots,s_k\\ b_1,\cdots,b_k \end{array} \big) =
 \int_0^1 \int_0^{u_1}\cdots \int_0^{u_{|\vec{s}|-1}} \frac{du_1}{g_1(u_1)}\cdots \frac{du_{|\vec{s}|}}{g_{|\vec{s}|}(u_{|\vec{s}|})},
 \mlabel{eq:intmp2}
\end{equation}
as commented in the introduction of~\mcite{3BL}.
Here
$$g_j(u_j)=\left\{\begin{array}{ll} b_i-u_j& \text{if } j= s_1+\cdots +s_i, 1\leq i \leq k,\\
u_j & \text{otherwise}. \end{array} \right .
$$

\noindent
The multiplication rule of two \MPVs according to their integral representations in Eq.~(\mref{eq:intmp2}) follows from the algebra homomorphism~\cite[\S 5.4]{3BL}
$$
 \Li\shf: \calh\shf\shzero(\ssg{S^1}) \to \mplalg,
 \quad x_0^{s_1-1}  x_{b_1}  \cdots   x_0^{s_k-1}  x_{b_k} \mapsto \lambda\big(\begin{array}{c} s_1,\cdots,s_k \\ b_1,\cdots, b_k\end{array}\big). $$

The algebra isomorphism $\rho:\calh\shf\lone(\ssg{S^1}) \to \calh\bsh(\zsg{S}^1)$ in Eq.~(\mref{eq:beta})
restricts to an algebra isomorphism
$$
\rho: \calh\shf\shzero(\ssg{S^1}) \to \calh\bsh\lizero (\zsg{S}^1),\quad
 x_0^{  \sind_1-1}  x_{b_1}  \cdots   x_0^{  \sind_k-1}  x_{b_k}
 \leftrightarrow \spair{ \sind_1,\  \sind_2,\ \cdots, \ \sind_k}{{b_1},\
 {b_2},\ \cdots,\ {b_k}}.
$$
Similarly the algebra isomorphisms $\eta:\calh\shf(\ssg{S^1}) \to
\calh\esh (\zsg{S}^1)$ in Eq.~(\mref{eq:shqsh1}) and $\theta: \calh\bsh(\zsg{S}^1) \to \calh\esh(\zsg{S}^1)$ in Eq.~(\mref{eq:la}) restrict to algebra isomorphisms
$$
\eta:\calh\shf\shzero(\ssg{S^1}) \to
\calh\esh\lizero (\zsg{S}^1), \quad
\theta: \calh\bsh\lizero(\zsg{S}^1) \to \calh\esh\lizero(\zsg{S}^1).
$$
Define
\begin{equation}
\Li\bsh: \calh\bsh\lizero(\zsg{S}^1) \to \mplalg, \quad
\spair{s_1,\cdots,s_k}{b_1,\cdots,b_k} \mapsto \lambda\big(\begin{array}{c} s_1,\cdots,s_k\\ b_1,\cdots,b_k \end{array} \big)
\mlabel{eq:lirho}
\end{equation}
and
\begin{equation}
\Li\esh: \calh\esh\lizero(\zsg{S}^1) \to \mplalg, \quad
\spair{s_1,\cdots,s_k}{z_1,\cdots,z_k} \mapsto \Li_{s_1,\cdots,s_k}(z_1,\cdots,z_k).
\mlabel{eq:lieta}
\end{equation}
Then we can organize these maps into the following commutative diagram extending the commutative diagram in~(\mref{eq:comm}):
$$
\xymatrix{
\calh\shf\shzero(\ssg{S^1}) \ar^{\rho}[rr]  \ar_{\Li\shf}[rrd] \ar@/^2pc/^{\eta}[rrrr]&&
    \calh\bsh\lizero(\zsg{S}^1) \ar^{\Li\bsh}[d]
    \ar^{\theta}[rr] && \calh\esh\lizero(\zsg{S}^1)
    \ar^{\Li\esh}[dll]
    \\
&& \mplalg &&
}
$$
where the commutativity of the left triangle follows from the definitions of the maps and that of the right triangle follows from Eq.~(\mref{eq:3b}).
Since $\Li\shf$ is an algebra homomorphism and $\rho$ and $\eta$ are algebra isomorphisms, it follows that $\Li\bsh$ and $\Li\esh$ are also algebra homomorphisms.

Therefore, applying $\Li\bsh$ to the two sides of Eq.~(\mref{eqn:maincoef}) in Theorem~\mref{thm:qshsh}, we obtain

\begin{coro}
Let $k,\ell$ be positive integers. Let $\vec{\rind}\in \ZZ_{\geq 1}^k$ and $\vec{\sind}\in \ZZ_{\geq
1}^\ell$. Let $\vec{a}=(a_1,\cdots,a_k)\in (S^1)^k$ and
$\vec{b}=(b_1,\cdots,b_\ell)\in (S^1)^\ell$ such that
$\spair{\rind_1}{a_1}\neq \spair{1}{1}$ and
$\spair{\sind_1}{b_1}\neq \spair{1}{1}$.
$$
\lambda\big(\begin{array}{c}\vec{\rind}\\ \vec{a}\end{array} \big)\, \lambda\big(\begin{array}{c}\vec{\sind}\\ \vec{b}\end{array} \big) =
\sum_{\vec{\tind}\in \ZZ_{\geq 1}^{k+\ell}, |\vec{\tind}|=
  |\vec{\rind}|+|\vec{\sind}|}
  \sum _{(\varphi,\psi)\in \indI_{k,\ell}}
  \Big(\prod_{i=1}^{k+\ell}c_{\vec{\rind},\vec{\sind}}^{\vec{\tind},
(\varphi,\psi)}(i)\Big)
    \lambda\big(\begin{array}{c} \vec{\tind}\\ \vec{a} \ssha_{(\varphi,\psi)} \vec{b}\end{array} \big).
$$
where $c_{\vec{\rind},\vec{\sind}}^{\vec{\tind}, (\varphi,\psi)}(i)$
is given in Eq.~(\mref{eq:coef-re-def1}) and
$\vec{a}\ssha_{(\varphi,\psi)}\vec{b}$ is given in
Eq.~(\mref{eq:mulind}). \mlabel{co:mainmpl1}
\end{coro}

Similarly, applying $\Li\esh$ to the two sides of Eq.~(\mref{eqn:maincoefe}) in Theorem~\mref{thm:qshshe}, we obtain

\begin{coro}
Let $k,\ell$ be positive integers. Let $\vec{\rind}\in \ZZ_{\geq 1}^k$ and $\vec{\sind}\in \ZZ_{\geq
1}^\ell$. Let $\vec{w}=(w_1,\cdots,w_k)\in (S^1)^k$ and
$\vec{z}=(z_1,\cdots,z_\ell)\in (S^1)^\ell$ such that
$\spair{\rind_1}{w_1}\neq \spair{1}{1}$ and
$\spair{\sind_1}{z_1}\neq \spair{1}{1}$. Then
$$
\Li_{\vec{\rind}}(\vec{w})\, \Li_{\vec{\sind}}(\vec{z}) =
\sum_{\vec{\tind}\in \ZZ_{\geq 1}^{k+\ell}, |\vec{\tind}|=
  |\vec{\rind}|+|\vec{\sind}|}
\,   \sum _{(\varphi,\psi)\in \indI_{k,\ell}}
  \Big(\prod_{i=1}^{k+\ell}c_{\vec{\rind},\vec{\sind}}^{\vec{\tind},
(\varphi,\psi)}(i)\Big)
    \Li_{\vec{\tind}}(\vec{w} \pssha_{(\varphi,\psi)} \vec{z})
$$
where $c_{\vec{\rind},\vec{\sind}}^{\vec{\tind}, (\varphi,\psi)}(i)$
is given in Eq.~(\mref{eq:coef-re-def1}) and
$\vec{w}\pssha_{(\varphi,\psi)}\vec{z}$ is given in
Eq.~(\mref{eq:pssha}). \mlabel{co:mainmpl}
\end{coro}
See Section~\mref{ss:exam} for examples in low dimensions.

\subsubsection{Multiple zeta values and alternating Euler sums}
\mlabel{sss:mzv}
Taking $z_i=1, 1\leq i\leq r,$ in
$
\Li_{s_1,\cdots,s_k}(z_1,\cdots,z_k)
$
defined in Eq.~(\mref{eq:mpls}) and the corresponding integral representation in Eq.~(\mref{eq:intmpl}), we obtain the \MZV and its integral representation:
\begin{eqnarray*}
 \zeta(s_1,\cdots,s_k):&= &\sum_{n_1>\cdots>n_k\geq 1} \frac{1}{n_1^{s_1}\cdots n_k^{s_k}}
 \\
 &=&
 \int_0^1 \int_0^{u_1}\cdots \int_0^{u_{|\vec{s}|-1}} \frac{du_1}{f_1(u_1)}\cdots \frac{du_{|\vec{s}|}}{f_{|\vec{s}|}(u_{|\vec{s}|})}
\end{eqnarray*}
for integers $s_i\geq 1$ and $s_1>1$. Here
$$f_j(u_j)=\left\{\begin{array}{ll} 1-u_j & \text{if } j= s_1,s_1+s_2,\cdots, s_1+\cdots +s_k,\\
u_j & \text{otherwise}. \end{array} \right .
$$
This is also the case when $\sg=\{1\}$ in our framework  in
Section~\mref{ss:frame} and \mref{ss:main}. Then we can identify
$\zsg{\sg}$ with $\ZZ_{\geq 1}$ and denote
$\vec{\nu}=\spair{s_1,\cdots,s_k}{z_1,\cdots,z_k}
=\spair{s_1,\cdots,s_k}{1,\cdots,1}$ by $\frakz=\frakz_{s_1}\cdots \frakz_{s_k}.$
Then
$\calh\qsh(\zsg{\sg})$ coincides with the quasi-shuffle algebra
$\calh\qsh$ encoding \MZVs~\mcite{Ho2,IKZ} through the identification $\frakz_{s_1}\cdots \frakz_{s_k}\leftrightarrow z_{s_1}\cdots z_{s_k}$. We will use $\frakz_{s_1}\cdots \frakz_{s_k}$ in place of $z_{s_1}\cdots z_{s_k}$ to avoid confusion with the vector $(z_1,\cdots,z_k)$ in $\vec{\nu}$.
$\calh\qsh$ contains the
subalgebra
$$
\calh\qsh\lzero:=\ZZ \oplus \ZZ\big\{\frakz_{s_1}\cdots \frakz_{s_k}\ \big|\
s_i\geq 1, s_1>1, 1\leq
i\leq k, k\geq 1 \big\}.
$$
Likewise the shuffle algebra $\calh\shf(\ssg{\sg})$, when $\sg=\{1\}$, coincides with the shuffle algebra $\calh\shf$~\mcite{Ho1,IKZ} encoding \MZVs, and there are subalgebras
$$
 \calh\shf\shzero:=\ZZ \oplus x_0\calh\shf x_1
\subseteq \calh\shf\lone:=\ZZ \oplus \calh\shf x_1
\subseteq \calh\shf,
$$
where $\calh\shf\lone$ coincides with $\calh\shf\lone(\zsg{\sg})$ defined in Eq.~(\mref{eq:shsubg}).
The natural isomorphism $\shqs: \calh\shf\lone \to \calh\qsh$ of abelian groups in Eq.~(\mref{eq:shqsh1})
restricts to an isomorphism of abelian groups
$$
\shqs: \calh\shf\shzero \to \calh\qsh\lzero, \quad
 x_0^{  s_1-1}  x_1  \cdots   x_0^{  s_k-1}  x_1
 \leftrightarrow \frakz_{s_1}\cdots \frakz_{s_k}.
$$
With the notation $ \frakz' \qssha \frakz'': = \shqs(
\shqs^{-1}(\frakz')\ssha \shqs^{-1}(\frakz'')) $ from
Eq.~(\mref{eq:shtrans}), the {\bf double shuffle relation} of \MZVs
is simply the ideal generated by the set
$$
\{ \frakz'\qssha \frakz'' - \frakz' \ast \frakz''\ |\
\frakz',\frakz''\in\calh\qsh\lzero\}$$ and the {\bf extended
double shuffle relation} of \MZVs~\mcite{IKZ} is the ideal generated
by the set
$$
\{ \frakz'\qssha \frakz'' - \frakz' \ast \frakz'',\ \frakz_1\qssha \frakz'' - \frakz_1 \ast \frakz''\ |\
\frakz',\frakz''\in\calh\qsh\lzero\}.$$

While the product $\frakz'\ast \frakz''$ simply follows from the
quasi-shuffle relation, the evaluation of  $\frakz'\qssha \frakz''$
involves first pulling $\frakz'$ and $\frakz''$ back to
$\calh\shf\shzero$ by $\shqs$, then expressing the shuffle product
$\shqs(\frakz') \ssha \shqs(\frakz'')$ as a linear combination of
words in $M(x_0,x_1)$, and then sending the result forward to
$\calh\qsh\lzero$ by $\shqs$. While this process can be defined
recursively (see Proposition~\mref{pp:rec}), the explicit formula
is found only in special cases, such as when $\frakz'=\frakz_{r},
\frakz''=\frakz_{s}$ are both of dimension one. As we have discussed in the
Introduction, the explicit formula in this case is Euler's formula
in Eq.~(\mref{eq:euler}).

Our Theorem~\mref{thm:qshshe} provides an explicit formula for $\qssha$ and hence for the shuffle product of \MZVs in the full generality.
\begin{coro}
Let $\vec{\rind}\in \ZZ_{\geq 1}^k$ and $\vec{\sind}\in \ZZ_{\geq
1}^\ell$ with $\rind_1,\sind_1\geq 2$. Then
$$
\zeta(\vec{\rind})\,\zeta(\vec{\sind}) =   \sum_{\vec{\tind}\in
\ZZ_{\geq 1}^{k+\ell}, |\vec{\tind}|=
  |\vec{\rind}|+|\vec{\sind}|}
  \Big(\sum _{(\varphi,\psi)\in \indI_{k,\ell}}
  \prod_{i=1}^{k+\ell}c_{\vec{\rind},\vec{\sind}}^{\vec{\tind},
(\varphi,\psi)}(i)\Big)
    \zeta(\vec{\tind})
$$
where $c_{\vec{\rind},\vec{\sind}}^{\vec{\tind}, (\varphi,\psi)}(i)$
is given in Eq.~(\mref{eq:coef-re-def1}).
\mlabel{co:mainmzv}
\end{coro}
See Section~\mref{ss:exam} for its specialization to Euler's decomposition formula and other special cases.
\begin{proof}
Since $\zeta(\vec{r})=\Li_{\vec{r}}(\vec{w})$ and $\zeta(\vec{s})=\Li_{\vec{s}}(\vec{z})$ where the vectors $\vec{w}$ and $\vec{z}$ have 1 as the components,
 the vectors $\vec{w}\pssha_{(\varphi,\psi)}
\vec{z}$ also have $1$ as their components and thus are independent of
the choice of $(\varphi,\psi)\in \indI_{k,\ell}$. Then the corollary
follows Corollary~\mref{co:mainmpl}.
\end{proof}

Between the case of \MZVs and the case of \MPVs, there is the case of {\bf alternating Euler sums}, defined by
$$
\zeta(s_1,\cdots,s_k;\sigma_1,\cdots,\sigma_k):=
\sum_{n_1>\cdots>n_k\geq 1} \frac{\sigma_1^{n_1}\cdots \sigma_k^{n_k}}{n_1^{s_1}\cdots n_k^{s_k}}\,,
$$
where $\sigma_i=\pm 1, 1\leq i\leq k.$
This corresponds to the case when $\sg=\{\pm 1\}$ in our framework.
More generally when $\sg$ is the group of $k$-th roots of unity, we have the {\bf multiple polylogarithms at roots of unity}~\mcite{Ra}.
We will not go into the details, but will give an example in Eq.~(\mref{eq:alteuler}) that generalizes Euler's formula.

\subsection{Examples}
\mlabel{ss:exam}
We now consider some special cases of Theorem~\mref{thm:qshshe}, Corollary~\mref{co:mainmpl} and Corollary~\mref{co:mainmzv}.
\subsubsection{The case of $r=s=1$}
In this case $\vec{\rind}=\rind_1$ and $\vec{\sind}=\sind_1$ are
positive integers, and $\vec{w}=w_1$ and $\vec{z}=z_1$ are in $\sg$.
Let $\vec{\tind}=(\tind_1, \tind_2)\in\ZZ_{\geq 1}^2$ with
$\tind_1+\tind_2=\rind_1+\sind_1$. If
$(\varphi,\psi)\in\indI_{1,1}$, then either $\varphi(1)=1$ and
$\psi(1)=2$, or $\psi(1)=1$ and $\varphi(1)=2$.  If $\varphi(1)=1$
and $\psi(1)=2$, then by Eq.~(\mref{eq:coef-re-def1}), we obtain
$$ c_{\rind_1,\sind_1}^{\vec{\tind}, (\varphi,\psi)}(1)
=
\binc{\tind_1-1}{\rind_1-1}, \quad
c_{\rind_1,\sind_1}^{\vec{\tind}
(\varphi,\psi)}(2)
=\binc{\tind_2-1}{
\tind_1+\tind_2-\rind_1-\sind_1}=1
$$ and thus
$$ c_{\rind_1,\sind_1}^{\vec{\tind}, (\varphi,\psi)}=
c_{\rind_1,\sind_1}^{\vec{\tind}, (\varphi,\psi)}(1)\,
c_{\rind_1,\sind_1}^{\vec{\tind}, (\varphi,\psi)}(2)
=\binc{\tind_1-1}{\rind_1-1}.$$ %Further, $\varphi^{-1}([1])=[1],
% \psi^{-1}([1])=\emptyset, \varphi^{-1}([2])=[1]$ and
% $\psi^{-1}([2])=[1]$. So
By Eq.~(\mref{eq:pssha}), we have
$$
\vec{w}\pssha_{(\varphi,\psi)} \vec{z} =(w_1,z_1/w_1).
$$

If $\psi(1)=1$ and $\varphi(1)=2$, then by Eq.~(\mref{eq:coef-re-def1}), we obtain
$$c_{\rind_1,\sind_1}^{\vec{\tind}, (\varphi,\psi)}(1)
= \binc{\tind_1-1}{\sind_1-1}, \quad
c_{\rind_1,\sind_1}^{\vec{\tind}, (\varphi,\psi)}(2)
=\binc{\tind_2-1}{
\tind_1+\tind_2-\rind_1-\sind_1}=1
$$ and thus
$$ c_{\rind_1,\sind_1}^{\vec{\tind}, (\varphi,\psi)}=
 c_{\rind_1,\sind_1}^{\vec{\tind}, (\varphi,\psi)}(1)\, c_{\rind_1,\sind_1}^{\vec{\tind}, (\varphi,\psi)}(2)=
\binc{\tind_1-1}{\sind_1-1}.$$ By Eq.~(\mref{eq:pssha}), we have $
\vec{w}\pssha_{(\varphi,\psi)} \vec{z} =(z_1,w_1/z_1). $
Therefore,
$$
\begin{aligned}
\spair{\rind_1}{w_1}\qssha \spair{\sind_1}{z_1} &= \hspace{-.5cm}
\sum_{\tind_1,\tind_2\geq 1, \tind_1+\tind_2=\rind_1+\sind_1}
\binc{\tind_1-1}{\rind_1-1} \spair{\tind_1,\tind_2}{w_1,z_1/w_1} +
\hspace{-.5cm}
\sum_{\tind_1,\tind_2\geq 1, \tind_1+\tind_2=\rind_1+\sind_1}
\binc{\tind_1-1}{\sind_1-1} \spair{\tind_1,\tind_2}{z_1,w_1/z_1} \\
&= \hspace{-.5cm}
\sum_{\tind_1,\tind_2\geq 1, \tind_1+\tind_2=\rind_1+\sind_1}
\binc{\tind_1-1}{t_1-r_1} \spair{\tind_1,\tind_2}{w_1,z_1/w_1} +
\hspace{-.5cm}
\sum_{\tind_1,\tind_2\geq 1, \tind_1+\tind_2=\rind_1+\sind_1}
\binc{\tind_1-1}{t_1-s_1} \spair{\tind_1,\tind_2}{z_1,w_1/z_1} \\
&=\sum_{k=0}^{\sind_1-1} \binc{\rind_1+k-1}{k}
\spair{\rind_1+k,\sind_1-k}{w_1,z_1/w_1} + \sum_{k=0}^{\rind_1-1}
\binc{\sind_1+k-1}{k} \spair{\sind_1+k,\rind_1-k}{z_1,w_1/z_1}
\end{aligned}
$$
by a change of variables $k=\tind_1-\rind_1$ for the first sum and
$k=\tind_1-\sind_1$ for the second sum.
Then by Corollary~\mref{co:mainmpl}, we obtain
the following relation for double polylogarithms
$$
{ \Li_{\rind_1}(w_1) \Li_{\sind_1}(z_1) \!
 = \! \sum_{k=0}^{\sind_1-1}\! \binc{\rind_1+k-1}{k} \Li_{\rind_1+k,\sind_1-k}(w_1,z_1/w_1)
+\! \sum_{k=0}^{\rind_1-1}\! \binc{\sind_1+k-1}{k}
\Li_{\sind_1+k,\rind_1-k}(z_1,w_1/z_1), }$$ where
$\rind_1,\sind_1\geq 1$, $w_1,z_1\in S^1$ and $(\rind_1,w_1)\neq
(1,1) \neq (\sind_1,z_1)$. In the special case when $w_1=\pm 1$ and
$z_1=\pm 1$, we have the following relation for alternating Euler
sums
\begin{equation}
\begin{aligned}
\zeta(\rind_1;w_1) \zeta(\sind_1;z_1)
 = &\sum_{k=0}^{s_1-1} \binc{\rind_1+k-1}{k} \zeta(\rind_1+k,\sind_1-k;w_1,z_1/w_1)\\
& + \sum_{k=0}^{\rind_1-1} \binc{\sind_1+k-1}{k}
\zeta(\sind_1+k,\rind_1-k;z_1,w_1/z_1),
\end{aligned}
\mlabel{eq:alteuler}
\end{equation}
when $\rind_1,\sind_1\geq 1$ and $(\rind_1,w_1)\neq (1,1) \neq
(\sind_1,z_1)$.

Further specializing, when $\rind_1,\sind_1\geq 2$ and $w_1=z_1=1$,
we obtain the decomposition formula of Euler in
Eq.~(\mref{eq:euler}).

\subsubsection{The case of $r=1,s=2$}
In this case $\spair{\vec{\rind}}{\vec{w}}=\spair{\rind_1}{w_1}$ and
$\spair{\vec{\sind}}{\vec{z}} =\spair{\sind_1,\sind_2}{z_1,z_2}$.
Let $\vec{\tind}=(\tind_1, \tind_2, \tind_3)\in\ZZ_{\geq 1}^3$ with
$\tind_1+\tind_2+\tind_3=\rind_1+\sind_1+\sind_2$. There are $3$
pairs $(\varphi,\psi)$ in $\indI_{1,2}$.

When $\varphi(1)=1$, $\psi(1)=2$ and $\psi(2)=3$, by Eq.~(\mref{eq:coef-re-def1}), we have
$$c_{\rind_1,\vec{\sind}}^{\vec{\tind}, (\varphi,\psi)}(1)
 = \binc{\tind_1-1}{\rind_1-1},
\quad
c_{\rind_1,\vec{\sind}}^{\vec{\tind}, (\varphi,\psi)}(2) =
\binc{\tind_2-1}{t_1+t_2-r_1-s_1}, \quad
c_{\rind_1,\vec{\sind}}^{\vec{\tind}, (\varphi,\psi)}(3) =
\binc{\tind_3-1}{ \sind_2-1}.
$$
When the second and the third terms are nonzero, we have $t_1+t_2\geq r_1+s_2$ and $t_3\geq s_2$. Then the inequalities must be equalities and we have
$c_{\rind_1,\vec{\sind}}^{\vec{\tind}, (\varphi,\psi)}(2)=
c_{\rind_1,\vec{\sind}}^{\vec{\tind}, (\varphi,\psi)}(3)=1.$
Thus
$$ c_{\rind_1,\vec{\sind}}^{\vec{\tind}, (\varphi,\psi)}
=c_{\rind_1,\vec{\sind}}^{\vec{\tind}, (\varphi,\psi)}(1)\,
c_{\rind_1,\vec{\sind}}^{\vec{\tind}, (\varphi,\psi)}(2)\,
c_{\rind_1,\vec{\sind}}^{\vec{\tind}, (\varphi,\psi)}(3)
=\left\{\begin{array}{cl} \binc{\tind_1-1}{\rind_1-1}, & \text{ if } t_3=s_2, \\ 0, & \text{ otherwise }. \end{array} \right .$$
By Eq.~(\mref{eq:pssha}) we have
$$
\vec{w}\pssha_{(\varphi,\psi)} \vec{z} = (w_1,z_1/w_1,z_2).
$$

Similarly, when $\varphi(1)=2$, $\psi(1)=1$ and $\psi(2)=3$, we have
$$ c_{\rind_1,\vec{\sind}}^{\vec{\tind}, (\varphi,\psi)}= \binc{\tind_1-1}{\sind_1-1}\binc{\tind_2-1}{s_2-t_3},\quad
\vec{w}\pssha_{(\varphi,\psi)} \vec{z} =(z_1,w_1/z_1,z_1z_2/w_1),
$$
and when $\varphi(1)=3$, $\psi(1)=1$ and $\psi(2)=2$, we have
$$ c_{\rind_1,\vec{\sind}}^{\vec{\tind}, (\varphi,\psi)}=
\binc{\tind_1-1}{\sind_1-1}\binc{\tind_2-1}{\sind_2-1}, \quad
\vec{w}\pssha_{(\varphi,\psi)} \vec{z} =(z_1,z_2,w_1/(z_1z_2)).
$$

Combining these computations with Corollary~\mref{co:mainmpl} we
obtain, for $\rind_1,\sind_1, \sind_2\geq 1$ and $(\rind_1,w_1)\neq
(1,1)\neq (\sind_1,z_1)$,
\begin{eqnarray*}
\Li_{\rind_1}(w_1)\,\Li_{\sind_1,\sind_2}(z_1,z_2)  &=
&
\hspace{-.5cm}\sum\limits_{
\substack{\tind_1,\tind_2,\tind_3\geq 1
\\ \tind_1+\tind_2
=\rind_1+\sind_1 }}   \!\!\!
\binc{\tind_1-1}{\rind_1-1} \Li_{(\tind_1,\tind_2,s_2)}
(w_1,z_1/w_1,z_2)
\\ && +
\sum_{
\substack{\tind_1,\tind_2,\tind_3\geq 1
\\ \tind_1+\tind_2+\tind_3 \\
=\rind_1+\sind_1+\sind_2 }}\bigg[
\binc{\tind_1-1}{\sind_1-1}\binc{\tind_2-1}{s_2-t_3}
\Li_{(\tind_1,\tind_2,\tind_3)} (z_1,w_1/z_1,z_1z_2/w_1)
\\
    &&  \qquad \qquad \qquad +
\binc{\tind_1-1}{\sind_1-1}\binc{\tind_2-1}{\sind_2-1}
\Li_{(\tind_1,\tind_2,\tind_3)} (z_1,z_2,w_1/(z_1z_2))\bigg ].
\end{eqnarray*}
Taking $w_1=z_1=z_2=1$ (or by Corollary~\mref{co:mainmzv}) we obtain the relation in Eq.~(\mref{eq:case1-2}) among \MZVs.

\subsubsection{The case of $r=s=2$}

In this case $\spair{\vec{\rind}}{\vec{w}}
=\spair{\rind_1,\rind_2}{w_1,w_2}$ and
$\spair{\vec{\sind}}{\vec{z}}=\spair{\sind_1, \sind_2}{z_1,z_2}$.
Let $\vec{\tind}=(\tind_1, \tind_2, \tind_3,\tind_4)\in\ZZ_{\geq
1}^4$ with
$\tind_1+\tind_2+\tind_3+\tind_4=\rind_1+\rind_2+\sind_1+\sind_2$.
Then there are $\binc{4}{2}=6$ choices of
$(\varphi,\psi)\in\indI_{2,2}$.

If $\varphi(1)=1$, $\varphi(2)=2$, $\psi(1)=3$ and $\psi(2)=4$, by Eq.~(\mref{eq:coef-re-def1}), we have
\begin{eqnarray*} &c_{\vec{\rind},\vec{\sind}}^{\vec{\tind}, (\varphi,\psi)}(1)
 = \binc{\tind_1-1}{\rind_1-1}, \quad
c_{\vec{\rind},\vec{\sind}}^{\vec{\tind}, (\varphi,\psi)}(2) =
\binc{\tind_2-1}{\rind_2-1},
\\
&c_{\vec{\rind},\vec{\sind}}^{\vec{\tind}, (\varphi,\psi)}(3) =
\binc{\tind_3-1}{
t_1+t_2+t_3-r_1-r_2-s_1}
=\binc{\tind_3-1}{\sind_2-\tind_4},
\quad
c_{\vec{\rind},\vec{\sind}}^{\vec{\tind}, (\varphi,\psi)}(4) =
\binc{\tind_4-1}{\sind_2-1}.
\end{eqnarray*}
When the third and fourth terms are nonzero, we have
$t_1+t_2+t_3\geq r_1+r_2+s_1$ and $t_4\geq s_2$. Hence they must be equalities and thus
$c_{\vec{\rind},\vec{\sind}}^{\vec{\tind}, (\varphi,\psi)}(3) =c_{\vec{\rind},\vec{\sind}}^{\vec{\tind}, (\varphi,\psi)}(4) =1.$
Then
$$ c_{\vec{\rind},\vec{\sind}}^{\vec{\tind}, (\varphi,\psi)}
 =c_{\vec{\rind},\vec{\sind}}^{\vec{\tind}, (\varphi,\psi)}(1)\,
c_{\vec{\rind},\vec{\sind}}^{\vec{\tind}, (\varphi,\psi)}(2)\,
c_{\vec{\rind},\vec{\sind}}^{\vec{\tind},
(\varphi,\psi)}(3)\,c_{\vec{\rind},\vec{\sind}}^{\vec{\tind},
(\varphi,\psi)}(4)
=\left\{\begin{array}{cl} \binc{\tind_1-1}{\rind_1-1}\binc{\tind_2-1}{\rind_2-1}, &
\text{ if } t_4=s_2, \\ 0, & \text{ otherwise}. \end{array} \right .
$$
Similarly, if
$\varphi(1)=3$, $\varphi(2)=4$, $\psi(1)=1$ and $\psi(2)=2$, then
$$c_{\vec{\rind},\vec{\sind}}^{\vec{\tind}, (\varphi,\psi)}
=\left\{\begin{array}{cl} \binc{\tind_1-1}{\sind_1-1}\binc{\tind_2-1}{\sind_2-1}, & \text{ if } t_4=r_2, \\ 0, & \text{ otherwise}. \end{array} \right .
$$
If $\varphi(1)=1$, $\varphi(2)=3$, $\psi(1)=2$ and $\psi(2)=4$, then
$$ c_{\vec{\rind},\vec{\sind}}^{\vec{\tind}, (\varphi,\psi)}
= \binc{\tind_1-1}{\rind_1-1}\binc{\tind_2-1}
{t_1+t_2-r_1-s_1}\binc{\tind_3-1}{\sind_2-\tind_4}.$$
If $\varphi(1)=2$, $\varphi(2)=4$, $\psi(1)=1$ and
$\psi(2)=3$, then
$$\begin{aligned} c_{\vec{\rind},\vec{\sind}}^{\vec{\tind}, (\varphi,\psi)}
= \binc{\tind_1-1}{ \sind_1-1}\binc{\tind_2-1}{
t_1+t_2-r_1-s_1}\binc{\tind_3-1}{\rind_2-\tind_4}.\end{aligned}$$
If $\varphi(1)=1$, $\varphi(2)=4$, $\psi(1)=2$ and $\psi(2)=3$, then
$$c_{\vec{\rind},\vec{\sind}}^{\vec{\tind}, (\varphi,\psi)} =
\binc{\tind_1-1}{r_1-1} \binc{\tind_2-1}{t_1+t_2-r_1-s_1}\binc{\tind_3-1}{\sind_2-1}.
$$
If $\varphi(1)=2$, $\varphi(2)=3$, $\psi(1)=1$ and
$\psi(2)=4$, then
$$ c_{\vec{\rind},\vec{\sind}}^{\vec{\tind}, (\varphi,\psi)}  =
\binc{\tind_1-1}{\sind_1-1} \binc{\tind_2-1}{t_1+t_2-r_1-s_1}\binc{\tind_3-1}{
\rind_2-1}.
$$
Then from Corollary~\mref{co:mainmzv}, we obtain Eq.~(\mref{eq:c2-2}). We likewise obtain formulas for the products
of double multiple polylogarithms and those of double alternating Euler sums.

\section{Preparational lemmas}
\mlabel{sec:lem}

In this section we prove some properties of the coefficients
$c^{\vec{\tind},(\varphi,\psi)}_{\vec{\rind},\vec{\sind}}$ in our
Theorem \mref{thm:qshsh} and Theorem~\mref{thm:qshshe} in preparation for their proofs in the next section.

We recall some notations from
Section \mref{ss:main}. Let $k,\ell\geq 1$, $\vec{\rind}\in
\ZZ_{\geq 1}^k$, $\vec{\sind}\in \ZZ_{\geq 1}^\ell$, $\vec{\tind}\in
\ZZ_{\geq 1}^{k+\ell}$ with
$|\vec{\tind}|=|\vec{\rind}|+|\vec{\sind}|$ and $(\varphi,\psi)\in
\indI_{k,\ell}$ be given. For $1\leq i\leq k+\ell$, denote
\begin{equation}
h_{(\varphi,\psi),i}=h_{(\varphi,\psi),(\vec{\rind},\vec{\sind}),i}=
       \left\{
              \begin{array}{ll} \rind_j & \text{ if } i=\varphi(j),
                                 \\
                                \sind_j & \text{ if } i=\psi(j).
              \end{array}
       \right.
\mlabel{eq:h}
\end{equation}
We note that, if we define
\begin{equation}
\vep_{\varphi,\psi}(i)=\left\{ \begin{array}{cl} 1 & \text{ if} \quad i\in
\im (\varphi), \\ -1 & \text{ if} \quad i\in \im (\psi),
\end{array}\right.
\mlabel{eq:vep1}
\end{equation}
then Eq.~(\mref{eq:coef-re-def1}) can be rewritten as
\begin{equation}
c_{\vec{\rind},\vec{\sind}}^{\vec{\tind},(\varphi,\psi)}(i) =\left\{
\begin{array}{ll}
\binc{\tind_i-1}{h_{(\varphi,\psi),i}-1} &
\begin{array}{l} \text{if } i=1 \\
\text{or if }
i\geq 2 \text{ and } \varepsilon_{\varphi,\psi}(i)\varepsilon_{\varphi,\psi}(i-1)=1,
\end{array}
\vspace{.2cm}
\\
\binc{\tind_{i}-1}{\sum\limits_{j=1}^{i} \tind_j
-\sum\limits_{j=1}^{i} h_{(\varphi,\psi),j}} & \text{ if } i\geq 2
\text{ and }
\varepsilon_{\varphi,\psi}(i)\varepsilon_{\varphi,\psi}(i-1)=-1.
\end{array}
\right.
\mlabel{eq:coef-re-def}
\end{equation}
Also recall
$$ c_{\vec{\rind},\vec{\sind}}^{\vec{\tind},(\varphi,\psi)}= \prod_{i=1}^{k+\ell}
c_{\vec{\rind},\vec{\sind}}^{\vec{\tind},(\varphi,\psi)}(i).$$

For the inductive proof to work, we also include the case when one
of $k$ or $\ell$ (but not both) is zero which corresponds to the
case when $\vec{\mu}$ or $\vec{\nu}\in \calh\qsh\lzero(\zsg{\sg})$
is the empty word $1$.
We will use the convention that $\ZZ_{\geq 1}^0=\{\bfe\}$ and denote
$|\bfe|=0$. When $k=0$, $\ell\geq 1$, we will also denote
$\vec{\rind}=\bfe$, denote $\bff:[k](=\emptyset)\to [k+\ell]=[\ell]$
and denote $\indI_{0,\ell}=\{(\bff,\id_{[\ell]})\}$.
Similarly, when $\ell=0$, $k\geq 1$, we denote $\vec{\sind}=\bfe$,
$\bff:[\ell]\to [k+\ell]=[k]$ and $\indI_{k,0}=\{(\id_{[k]},
\bff)\}$. Then the notations in Eq.~(\mref{eq:h}) --
(\mref{eq:coef-re-def}) still make sense even if exactly one of $k$
and $\ell$ is zero. More precisely, when $k=0$, $\ell\geq 1$, we have
$h_{(\bff, \id_{[\ell]}),(\bfe, \vec{\sind}), i}=\sind_i$,
$\vep_{\bff,\id_{[\ell]}}(i)=-1, 1\leq i\leq \ell$. Also, for any
$\vec{\sind}$ and $\vec{\tind}\in \ZZ_{\geq 1}^\ell$ with
$|\vec{\sind}|=|\vec{\tind}|$, we have
\begin{equation}
c_{\bfe,\vec{\sind}}^{\vec{\tind},(\bff,\id_{[\ell]})}=
\prod_{i=1}^\ell\binc{\tind_i-1}{\sind_i-1} % =\delta_{\vec{\sind}}^{\vec{\tind}}
=\prod_{i=1}^\ell\delta_{\sind_i}^{\tind_i}.
\mlabel{eq:memp}\end{equation} Similarly, if $\vec{\sind}=\bfe$,
then for any $\vec{\rind},\vec{\tind}\in\ZZ_{\geq 1}^k$ with
$|\vec{\rind}|=|\vec{\tind}|$, we have
$h_{(\id_{[k]}, \bff),(\vec{r}, \bfe), i}=r_i$,
$\vep_{\id_{[k]},\bff}(i)=1, 1\leq i\leq k$ and
\begin{equation}
c_{\vec{\rind},\bfe}^{\vec{\tind},(\id_{[k]},\bff)}=
\prod_{i=1}^k\delta_{\rind_i}^{\tind_i}.
\mlabel{eq:nemp}\end{equation}

We first give some conditions for the vanishing of
$c^{\vec{\tind},(\varphi,\psi)}_{\vec{\rind},\vec{\sind}}$.
\begin{lemma} Let $k,\ell\geq 1$.
Let $\vec{\rind}\in \ZZ_{\geq 1}^k$, $\vec{\sind}\in \ZZ_{\geq
1}^\ell$ and $\vec{\tind}\in \ZZ_{\geq 1}^{k+\ell}$ with
$|\vec{\rind}|+|\vec{\sind}|=|\vec{\tind}|$. Let
$(\varphi,\psi)\in\indI_{k,\ell}$.
 Then
$c^{\vec{\tind},(\varphi,\psi)}_{\vec{\rind},\vec{\sind}}\neq 0$ if
and only if, for $1\leq i\leq k+\ell$,
$$ \left \{
\begin{array}{ll}
\tind_i\geq h_{(\varphi,\psi),i}, & \text{ if } i=1 \text{ or if }
i\geq 2 \text{ and } \varepsilon_{\varphi,\psi}(i)\varepsilon_{\varphi,\psi}(i-1)=1,\\
\sum\limits_{j=1}^{i}\tind_j \geq \sum\limits_{j=1}^{i}
h_{(\varphi,\psi),j} > \sum\limits_{j=1}^{i-1} \tind_j, & \text{ if
} i\geq 2 \text{ and }
\varepsilon_{\varphi,\psi}(i)\varepsilon_{\varphi,\psi}(i-1)=-1.
\end{array}
\right .$$
\mlabel{lem:coefzero}
\end{lemma}
\begin{proof}
By definition,
$c^{\vec{\tind},(\varphi,\psi)}_{\vec{\rind},\vec{\sind}}\neq 0$ if
and only if
$c^{\vec{\tind},(\varphi,\psi)}_{\vec{\rind},\vec{\sind}}(i)\neq 0$
for every $i\in [k+\ell]$. Also $\binc{a}{b}\neq 0$ if and only if
$a\geq b\geq 0$. Then the lemma follows since
$$\Big(\tind_i-1\geq
h_{(\varphi,\psi),i}-1\geq 0\Big) \Leftrightarrow \Big(\tind_i\geq
h_{(\varphi,\psi),i}\geq 1 \Big) $$ and
$$\Big(\tind_i-1\geq \sum_{j=1}^i\tind_j-\sum_{j=1}^ih_{(\varphi,\psi),i} \geq 0\Big)
\Leftrightarrow \Big(-\sum_{j=1}^{i-1} \tind_j > -\sum_{j=1}^{i-1}
\tind_j-1 \geq -\sum_{j=1}^ih_{(\varphi,\psi),j} \geq
-\sum_{j=1}^i\tind_j\Big).
$$
\end{proof}

\begin{lemma} Let $k,\ell, \vec{\rind},\vec{\sind},\vec{\tind}$ be as in
Lemma \mref{lem:coefzero}.
\begin{enumerate}
\item
Let $(\varphi,\psi)\in\indI_{k,\ell}$. If $\varphi(1)=1$,
$\sind_1=1$ and $\tind_1>\rind_1$ or if $\psi(1)=1$, $\rind_1=1$ and
$\tind_1>\sind_1$, then
$$c^{\vec{\tind}, (\varphi,\psi)}_{\vec{\rind},\vec{\sind}}=0. $$
\mlabel{it:coefzero1}
\item
If $\tind_1< \min (\rind_1, \sind_1)$, then
$c^{\vec{\tind},(\varphi,\psi)}_{\vec{\rind},\vec{\sind}}=0$ for any
$(\varphi,\psi)\in\indI_{k,\ell}$. \mlabel{it:coefzero2}
\end{enumerate}
\mlabel{lem:coefzero2}
\end{lemma}
\begin{proof}
(\mref{it:coefzero1}). We only consider the case when $\varphi(1)=1,
\sind_1=1$ and $\tind_1>\rind_1$. The proof of the other case is
similar. Since $\varphi(1)=1$, we have $\psi(1)>1$. This means that
$h_{(\varphi,\psi),i}=\rind_i$ for $1\leq i\leq \psi(1)-1$ and
$h_{(\varphi,\psi),\psi(1)}=\sind_1$. Suppose
$c^{\vec{\tind},(\varphi,\psi)}_{\vec{\rind},\vec{\sind}}\neq 0$.
Then by Lemma~\mref{lem:coefzero}, we have $\tind_i \geq \rind_i$
for $ 2\leq i \leq \psi(1)-1$ and $\sum\limits_{j=1}^{\psi(1)-1}
\rind_j +\sind_1
> \sum\limits_{j=1}^{\psi(1)-1} \tind_j$ by taking $i=\psi(1)$. From
these two inequalities, we obtain $\rind_1+\sind_1>\tind_1$ and
hence $\rind_1\geq \tind_1$ since $\sind_1=1$. This is a
contradiction.
\smallskip

\noindent (\mref{it:coefzero2}) If $\tind_1<\min(\rind_1,\sind_1)$,
then $\tind_1<h_{(\varphi,\psi),1}$. So by Lemma
\mref{lem:coefzero}, for every $(\varphi,\psi)\in \indI_{k,\ell}$ we
have $c_{\vec{\rind},\vec{\sind}}^{\vec{\tind},(\varphi,\psi)}=0$.
\end{proof}

We next give some relations among the numbers
$c^{\vec{\tind},(\varphi,\psi)}_{\vec{\rind},\vec{\sind}}(i)$ as the
parameters vary.

\begin{defn}
Let $\vec{e}_1$ denote $(1,0,\cdots,0)$ of suitable dimension. So
for any vector $\vec{x}=(x_1,x_2,\cdots,x_k)$ and $a\in \ZZ$, we have
$$\vec{x}-a
\vec{e}_1=(x_1-a,x_2,\cdots,x_k).$$
Define
$$\vec{x}\,'=(x_1',\cdots,x_{k-1}'):
=(x_2,\cdots,x_k)$$ with the convention that $(x_1)'=\bfe.$ For a
function $f$ on $[k]$, let $\dt{f}$ and $\st{f}$ be respectively the
functions on $[k-1]$ and $[k]$ defined by
$$
\dt{f}(x)=f(x+1)-1, \quad \st{f}(x)=f(x)-1 $$ with the convention
that $[0]=\emptyset$ and that, if $f$ is a function on $[1]$, then
$\dt{f}=\bff$. Let $f^{\&}$ and $f^*$ be respectively the functions
on $[k+1]$ and $[k]$ defined by
$$
f^{\&}(1)=1, \; f^{\&}(x)=f(x-1)+1,  \quad f^*(y)=f(y)+1, \quad
2\leq x\leq r+1, 1\leq y\leq r.
$$
\mlabel{de:notn}
\end{defn}

Also define
$$\indI_{k,\ell, \varphi(1)=1}=\{(\varphi,\psi)\in\indI_{k,\ell}\ |\
\varphi(1)=1\}, \quad \indI_{k,\ell,
\psi(1)=1}=\{(\varphi,\psi)\in\indI_{k,\ell}\ |\ \psi(1)=1\}.$$

\begin{lemma} Let $k,\ell\geq 1$. The map
$$(\dt{}, \st{}): \indI_{k,\ell,\varphi(1)=1}\rightarrow
\indI_{k-1,\ell}, \quad (\varphi,\psi) \mapsto
(\dt{\varphi},\st{\psi})$$ is a bijection whose inverse is given by
$$ ({^{\&}}, {^*}): \indI_{k-1,\ell}\rightarrow
\indI_{k,\ell,\varphi(1)=1}, \quad (\varphi,\psi)\mapsto
(\varphi^{\&},\psi^*).$$ Similarly, the map
$$(\st{}, \dt{}): \indI_{k,\ell,\psi(1)=1}\rightarrow
\indI_{k,\ell-1}, \quad (\varphi,\psi)\mapsto
(\st{\varphi},\dt{\psi})$$ is a bijection whose inverse is given by
$$ ({{^*},^{\&}}): \indI_{k,\ell-1}\rightarrow
\indI_{k,\ell,\psi(1)=1}, \quad (\varphi,\psi)\mapsto
(\varphi^{*},\psi^{\&}).$$ \mlabel{lem:bij}
\end{lemma}
\begin{proof} From
the definition we verify that $$({\dt{}}, {\st{}})(
\indI_{k,\ell,\varphi(1)=1})\subseteq \indI_{k-1,\ell}$$ and
$$ ({^{\&}}, {^*})(\indI_{k-1,\ell})\subseteq
\indI_{k,\ell,\varphi(1)=1}.$$ Then to prove the first assertion we
only need to show that $(\dt{\varphi})^{\&}=\varphi$ and
$(\st{\psi})^{*}=\psi$ if $\varphi(1)=1$, and that
$\dt{(\varphi^{\&})}=\varphi$ and $\st{(\psi^*)}=\psi$. We just
check the first equation and leave the others to the interested
reader. First we have $(\dt{\varphi})^{\&}(1)=1 $ by definition.
Since $\varphi(1)=1$, we have $(\dt{\varphi})^{\&}(i)=\varphi(i)$
when $i=1$. If $i\geq 2$, then by definition we have
$\dt{\varphi}(i-1)=\varphi(i)-1$ and
$(\dt{\varphi})^{\&}(i)=\dt{\varphi}(i-1)+1=\varphi(i)$, as desired.

The proof of the second assertion in the lemma is similar.
\end{proof}

\begin{lemma} Let $k,\ell, \vec{\rind},\vec{\sind},\vec{\tind}$ and $(\varphi,\psi)$ be as in
Lemma \mref{lem:coefzero}.
\begin{enumerate}
\item
Let $a$ and $b$ be integers such that $a<\min(\tind_1, \rind_1)$,
$b<\min (\tind_1, \sind_1)$. Then  for all $i\in\{2,\cdots,
k+\ell\}$, we have
$$ c_{\vec{\rind}-a\vec{e}_1, \vec{\sind}}^{\vec{\tind}-a\vec{e}_1,
(\varphi,\psi)}(i) =c_{\vec{\rind}, \vec{\sind}}^{\vec{\tind},
(\varphi,\psi)}(i)
$$  and
$$ c_{\vec{\rind}, \vec{\sind}-b\vec{e}_1}^{\vec{\tind}-b\vec{e}_1,
(\varphi,\psi)}(i) =c_{\vec{\rind}, \vec{\sind}}^{\vec{\tind},
(\varphi,\psi)}(i).
$$
\mlabel{it:coef1a}
\item
If $\varphi(1)=1$ and $\rind_1=\tind_1=1$, then
\begin{equation}
c_{\vec{\rind},\vec{\sind}}^{\vec{\tind},(\varphi,\psi)}(i+1)=
c_{\vec{\rind}\,',\vec{\sind}}^{\vec{\tind}\,',(\dt{\varphi},\st{\psi})}(i),\quad
1\leq i\leq k+\ell-1, \mlabel{eq:c1b}
\end{equation}
with the notations in Definition~\mref{de:notn}.
Similarly, if $\psi(1)=1$ and $\sind_1=\tind_1=1$, then
\begin{equation}
c_{\vec{\rind},\vec{\sind}}^{\vec{\tind},(\varphi,\psi)}(i+1)=
c_{\vec{\rind},\vec{\sind}\,'}^{\vec{\tind}\,',(\st{\varphi},\dt{\psi})}(i),\quad
1\leq i\leq k+\ell-1. \mlabel{eq:c2b}
\end{equation}
\mlabel{it:coef1b}
\end{enumerate}
\mlabel{lem:coefind}
\end{lemma}
\begin{proof}
(\mref{it:coef1a}) We prove the first equality. The proof for the
second equality is similar. Since $a<\min(\rind_1,\tind_1)$, we have
$\vec{\rind}-a\vec{e}_1\in\ZZ_{\geq 1}^k$ and
$\vec{\tind}-a\vec{e}_1\in \ZZ_{\geq 1}^{k+\ell}$. For better
distinction, we will use the full notation
$h_{(\varphi,\psi),(\vec{\rind},\vec{\sind}),i}$ defined in
Eq.~(\mref{eq:h}) instead of its abbreviation
$h_{(\varphi,\psi),i}$. Then we have
\begin{equation}
h_{(\varphi,\psi),(\vec{\rind}-a\vec{e}_1,\vec{\sind}),i}
=\left\{\begin{array}{ll}
h_{(\varphi,\psi),(\vec{\rind},\vec{\sind}),i} & \text{ if } i\neq
\varphi(1),
\\ h_{(\varphi,\psi),(\vec{\rind},\vec{\sind}),i}-a & \text{ if } i=\varphi(1).\end{array}\right.
\mlabel{eq:hs}
\end{equation}

Let $i\in\{ 2,\cdots, k+\ell\}$. If
$\vep_{\varphi,\psi}(i)\vep_{\varphi,\psi}(i-1)=1$, then $i\neq
\varphi(1)$. Indeed, if $i=\varphi(1)$, then $i-1$ must be in
$\im(\psi)$, implying that
$\vep_{\varphi,\psi}(i)\vep_{\varphi,\psi}(i-1)=-1$. Thus
$$c_{\vec{\rind}-a\vec{e}_1,\vec{\sind}}^{\vec{\tind}-a\vec{e}_1,(\varphi,\psi)}(i)= \binc{\tind_i-1}{h_{(\varphi,\psi),(\vec{\rind}-a\vec{e}_1,\vec{\sind}),i}-1}
=\binc{\tind_i-1}{h_{(\varphi,\psi),(\vec{\rind},\vec{\sind}),i}-1}
=c_{\vec{\rind},\vec{\sind}}^{\vec{\tind},(\varphi,\psi)}(i).$$ If
$\vep_{\varphi,\psi}(i)\vep_{\varphi,\psi}(i-1)=-1$, then either
$i=\varphi(j)$ or $i-1=\varphi(j)$ for some $j\in [k]$. In either
case, we have $i\geq \varphi(1)$ since $\varphi$ keeps the order.
Thus by Eq.~(\mref{eq:hs}), we have
$$\sum_{j=1}^i h_{(\varphi,\psi),(\vec{\rind}-a\vec{e}_1,\vec{\sind}),j}
=\sum_{j=1}^i h_{(\varphi,\psi),(\vec{\rind},\vec{\sind}),j}-a.$$ So
$$c_{\vec{\rind}-a\vec{e}_1,\vec{\sind}}^{\vec{\tind}-a\vec{e}_1,(\varphi,\psi)}(i) =
\binc{\tind_i-1}{(\tind_1-a)+\sum\limits_{j=2}^i\tind_j-\sum\limits_{j=1}^i
h_{(\varphi,\psi),(\vec{\rind}-a\vec{e}_1,\vec{\sind}),j}}
=\binc{\tind_i-1}{\sum\limits_{j=1}^i\tind_j-\sum\limits_{j=1}^i
h_{(\varphi,\psi),(\vec{\rind},\vec{\sind}),j}}=c_{\vec{\rind},\vec{\sind}}^{\vec{\tind},(\varphi,\psi)}(i).$$
\smallskip

\noindent (\mref{it:coef1b}) Let $\varphi(1)=1$ and
$\rind_1=\tind_1=1$. By Eq.~(\mref{eq:h}), for $1\leq i\leq
k+\ell-1$,
$$
\begin{aligned}
h_{(\varphi,\psi),(\vec{\rind},\vec{\sind}),i+1} &=
\left \{\begin{array}{ll} \rind_j & \text{ if } i+1=\varphi(j) \\
\sind_j & \text{ if } i+1=\psi(j) \end{array} \right .
= \left\{\begin{array}{ll} r'_{j-1} & \text{ if } i=\varphi(j)-1 \\
\sind_j & \text{ if } i=\psi(j)-1 \end{array} \right . \\
&= \left\{\begin{array}{ll} r'_{j} & \text{ if } i=\varphi(j+1)-1 \\
\sind_j & \text{ if } i=\psi(j)-1 \end{array} \right .
= \left\{\begin{array}{ll} r'_{j} & \text{ if } i=\dt{\varphi}(j) \\
\sind_j & \text{ if } i=\st{\psi}(j). \end{array} \right .
\end{aligned}
$$
Thus
\begin{equation}
h_{(\varphi,\psi),(\vec{\rind},\vec{\sind}),i+1} =
h_{(\dt{\varphi},\st{\psi}),(\vec{\rind}\,',\vec{\sind}),i}, 1\leq
i\leq k+\ell-1. \mlabel{eq:ch}
\end{equation}
Also, for $1\leq i\leq k+\ell-1$, since $\varphi(1)=1$, we have
$$
i+1\in \im(\varphi) \Leftrightarrow i+1=\varphi(j), j\in
\{2,\cdots,k\} \Leftrightarrow i=\dt{\varphi}(j-1), j-1\in [k-1]
\Leftrightarrow i\in \im (\dt{\varphi}).$$ Similarly, $i+1\in
\im(\psi) \Leftrightarrow i\in \im(\st{\psi}).$ Thus
\begin{equation}
\vep_{\varphi,\psi}(i+1)=\vep_{\dt{\varphi},\st{\psi}}(i), 1\leq
i\leq k+\ell-1. \mlabel{eq:ce}
\end{equation}

We now verify Eq.~(\mref{eq:c1b}) for $i=1$. Since $\varphi(1)=1$, either $2=\varphi(2)$ or $2=\psi(1)$.
If $2=\varphi(2)$, then
$\vep_{\varphi,\psi}(2)\vep_{\varphi,\psi}(1)=1$ and so
$$
c_{\vec{\rind},\vec{\sind}}^{\vec{\tind},(\varphi,\psi)}(2) =
\binc{\tind_2-1}{\rind_2-1} = \binc{\tind'_1-1}{r'_1-1}
=c_{\vec{\rind}\,',\vec{\sind}}^{\vec{\tind}\,',(\dt{\varphi},\st{\psi})}(1).
$$
If $\psi(1)=2$, then
$\vep_{\varphi,\psi}(2)\vep_{\varphi,\psi}(1)=-1$. So by the
condition that $\rind_1=\tind_1=1$, we obtain
$$c_{\vec{\rind},\vec{\sind}}^{\vec{\tind},(\varphi,\psi)}(2)
=\binc{\tind_2-1}{\tind_1+\tind_2-\rind_1-\sind_1} =
\binc{\tind_2-1}{\tind_2-\sind_1}=\binc{\tind_2-1}{\sind_1-1}
=\binc{\tind'_1-1}{\sind_1-1} =
c_{\vec{\rind}\,',\vec{\sind}}^{\vec{\tind}\,',(\dt{\varphi},\st{\psi})}(1).$$

Next consider $i\geq 2$. By Eq.~(\mref{eq:ch}) and Eq.~(\mref{eq:ce}), we have
$$
\begin{aligned}
c_{\vec{\rind},\vec{\sind}}^{\vec{\tind},(\varphi,\psi)}(i+1)
&=\left\{
\begin{array}{ll}
\binc{\tind_{i+1}-1}{h_{(\varphi,\psi),(\vec{\rind},\vec{\sind}),i+1}-1}
& \text{ if }
 \varepsilon_{\varphi,\psi}(i+1)\varepsilon_{\varphi,\psi}(i)=1,
\vspace{.2cm}
\\
\binc{\tind_{i+1}-1}{\sum\limits_{j=1}^{i+1} \tind_j
-\sum\limits_{j=1}^{i+1}
h_{(\varphi,\psi),(\vec{\rind},\vec{\sind}),j}} & \text{ if }
\varepsilon_{\varphi,\psi}(i+1)\varepsilon_{\varphi,\psi}(i)=-1,
\end{array}
\right. \\
&=\left\{
\begin{array}{ll}
\binc{\tind'_i-1}{h_{(\dt{\varphi},\st{\psi}),(\vec{\rind}\,',\vec{\sind}),
i}-1} & \text{ if }  \varepsilon_{\dt{\varphi},\st{\psi}}(i)
\varepsilon_{\dt{\varphi},\st{\psi}}(i-1)=1, \vspace{.2cm}
\\
\binc{\tind'_{i}-1}{\sum\limits_{j=1}^{i} \tind'_j
-\sum\limits_{j=1}^{i}
h_{(\dt{\varphi},\st{\psi}),(\vec{\rind}\,',\vec{\sind}),j}} & \text{
if } \varepsilon_{\dt{\varphi},\st{\psi}}(i)
\varepsilon_{\dt{\varphi},\st{\psi}}(i-1)=-1,
\end{array}
\right.
\end{aligned}
$$
since $\tind_1=1$ and
$h_{(\varphi,\psi),(\vec{\rind},\vec{\sind}),1}=\rind_1=1$.
Therefore, we have
$c_{\vec{\rind},\vec{\sind}}^{\vec{\tind},(\varphi,\psi)}(i+1)
=c_{\vec{\rind}\,',\vec{\sind}}^{\vec{\tind}\,',(\dt{\varphi},\st{\psi})}(i)$
when $i\geq 2$.

The proof for Eq.~(\ref{eq:c2b}) is similar.
\end{proof}

\begin{lemma} Let $k,\ell, \vec{\rind},\vec{\sind},\vec{\tind}$ and $(\varphi,\psi)$ be as in
Lemma \mref{lem:coefzero}.
\begin{enumerate}
\item Suppose that $\rind_1\geq 2$ and $\sind_1\geq 2$. If $\tind_1\geq 2$,
then we have
\begin{equation}
c_{\vec{\rind},\vec{\sind}}^{\vec{\tind},(\varphi,\psi)} =
c_{\vec{\rind}-\vec{e}_1,\vec{\sind}}^{\vec{\tind}-\vec{e}_1,(\varphi,\psi)}
+
c_{\vec{\rind},\vec{\sind}-\vec{e}_1}^{\vec{\tind}-\vec{e}_1,(\varphi,\psi)}
\mlabel{eq:coef1}
\end{equation}
If $\tind_1=1$, then we have
\begin{equation}
c_{\vec{\rind},\vec{\sind}}^{\vec{\tind},(\varphi,\psi)}=0.
\mlabel{eq:coefzero1}
\end{equation}
\mlabel{it:coef1}
\item
Suppose that $\rind_1=\sind_1=1$. If $\varphi(1)=1$ and $\tind_1=1$,
then we have
\begin{equation}
c_{\vec{\rind},\vec{\sind}}^{\vec{\tind},
(\varphi,\psi)}=c_{\vec{\rind}\,',\vec{\sind}}^{\vec{\tind}\,',(\dt{\varphi},
\st{\psi})} \mlabel{eq:coef21}
\end{equation} with the notations in Definition~\mref{de:notn}.
If $\psi(1)=1$ and $\tind_1=1$, then we have
\begin{equation}
c_{\vec{\rind},\vec{\sind}}^{\vec{\tind},
(\varphi,\psi)}=c_{\vec{\rind},\vec{\sind}\,'}^{\vec{\tind}\,',(\st{\varphi},
\dt{\psi})}. \mlabel{eq:coef22}
\end{equation}
If $\tind_1\geq 2$, then we have
\begin{equation}
c_{\vec{\rind},\vec{\sind}}^{\vec{\tind}, (\varphi,\psi)}=0.
\mlabel{eq:coefzero2}
\end{equation}
\mlabel{it:coef2}
\item
Suppose that $\rind_1=1$ and $\sind_1\geq 2$. If $\varphi(1)=1$ and
$\tind_1=1$, then we have
\begin{equation} c_{\vec{\rind},\vec{\sind}}^{\vec{\tind},
(\varphi,\psi)}=c_{\vec{\rind}\,',\vec{\sind}}^{\vec{\tind}\,',(\dt{\varphi},
\st{\psi})}. \mlabel{eq:coef31}
\end{equation}
If $\psi(1)=1$ and $\tind_1=1$, then we have
\begin{equation}
c_{\vec{\rind},\vec{\sind}}^{\vec{\tind}, (\varphi,\psi)}=0.
\mlabel{eq:coefzero3}
\end{equation} If $\tind_2\geq 2$, then we have
\begin{equation}
c_{\vec{\rind},\vec{\sind}}^{\vec{\tind},
(\varphi,\psi)}=c_{\vec{\rind},\vec{\sind}-\vec{e}_1}^{\vec{\tind}-\vec{e}_1,
(\varphi,\psi)}. \mlabel{eq:coef32}
\end{equation}
Similar statements hold when $\rind_1\geq 2$ and $\sind_1=1$.
\mlabel{it:coef3}
\end{enumerate}
\mlabel{lem:coef}
\end{lemma}
\begin{proof}
(\mref{it:coef1}) If $\varphi(1)=1$, then
$$c^{\vec{\tind},(\varphi,\psi)}_{\vec{\rind},\vec{\sind}}(1)=
\binc{\tind_1-1}{\rind_1-1}= \binc{\tind_1-2}{\rind_1-2}
+\binc{\tind_1-2}{\rind_1-1}
=c^{\vec{\tind}-\vec{e}_1,(\varphi,\psi)}_{\vec{\rind}-\vec{e}_1,\vec{\sind}}(1)
+c^{\vec{\tind}-\vec{e}_1,(\varphi,\psi)}_{\vec{\rind},\vec{\sind}-\vec{e}_1}(1).
$$ Similarly, if $\psi(1)=1$, we also have
$$c^{\vec{\tind},(\varphi,\psi)}_{\vec{\rind},\vec{\sind}}(1)
=c^{\vec{\tind}-\vec{e}_1,(\varphi,\psi)}_{\vec{\rind}-\vec{e}_1,\vec{\sind}}(1)
+c^{\vec{\tind}-\vec{e}_1,(\varphi,\psi)}_{\vec{\rind},\vec{\sind}-\vec{e}_1}(1).$$
In either case, by Lemma \mref{lem:coefind}.(\mref{it:coef1a}) we
have
$$c^{\vec{\tind},(\varphi,\psi)}_{\vec{\rind},\vec{\sind}}(i)
=c^{\vec{\tind}-\vec{e}_1,(\varphi,\psi)}_{\vec{\rind}-\vec{e}_1,\vec{\sind}}(i)
=c^{\vec{\tind}-\vec{e}_1,(\varphi,\psi)}_{\vec{\rind},\vec{\sind}-\vec{e}_1}(i)$$
when $i\in \{ 2,\cdots, k+\ell \}$. Hence {\allowdisplaybreaks
$$
\begin{aligned}
c^{\vec{\tind},(\varphi,\psi)}_{\vec{\rind},\vec{\sind}} & =
\prod_{i=1}^{k+\ell}c^{\vec{\tind},(\varphi,\psi)}_{\vec{\rind},\vec{\sind}}(i)
\\
&=
(c^{\vec{\tind}-\vec{e}_1,(\varphi,\psi)}_{\vec{\rind}-\vec{e}_1,\vec{\sind}}(1)
+c^{\vec{\tind}-\vec{e}_1,(\varphi,\psi)}_{\vec{\rind},\vec{\sind}-\vec{e}_1}(1))
\prod_{i=2}^{k+\ell}c^{\vec{\tind},(\varphi,\psi)}_{\vec{\rind},\vec{\sind}}(i)
\\
&= \prod_{i=1}^{k+\ell}
c^{\vec{\tind}-\vec{e}_1,(\varphi,\psi)}_{\vec{\rind}-\vec{e}_1,\vec{\sind}}(i)
+ \prod_{i=1}^{k+\ell}
c^{\vec{\tind}-\vec{e}_1,(\varphi,\psi)}_{\vec{\rind},\vec{\sind}-\vec{e}_1}(i)
\\
& = c^{\vec{\tind}-\vec{e}_1,
   (\varphi,\psi)}_{\vec{\rind}-\vec{e}_1,\vec{\sind}}+c^{\vec{\tind}-\vec{e}_1,
   (\varphi,\psi) }_{\vec{\rind},\vec{\sind}-\vec{e}_1}.
\end{aligned}
$$} This proves Eq.~(\mref{eq:coef1}). Eq~(\mref{eq:coefzero1})
follows from Lemma~\mref{lem:coefzero2} (\mref{it:coefzero2}).

\smallskip

\noindent (\mref{it:coef2}) First we assume that $\tind_1=1$. For
$(\varphi,\psi)\in\indI_{k,\ell}$, either $\varphi(1)=1$ or
$\psi(1)=1$. If $\varphi(1)=1$, then
$$c^{\vec{\tind},
(\varphi,
\psi)}_{\vec{\rind},\vec{\sind}}(1)=\binc{\tind_1-1}{\rind_1-1}=\binc{0}{0}=1$$
and by Lemma \mref{lem:coefind}.(\mref{it:coef1b}) we have
$$c^{\vec{\tind}, (\varphi,
\psi)}_{\vec{\rind},\vec{\sind}}(i+1)=
c_{\vec{\rind}\,',\vec{\sind}}^{\vec{\tind}\,',(\dt{\varphi},\st{\psi})}(i).
$$
Hence
$$
c_{\vec{\rind},\vec{\sind}}^{\vec{\tind},(\varphi,\psi)}
=\prod_{i=1}^{k+\ell}c_{\vec{\rind},\vec{\sind}}^{\vec{\tind},(\varphi,\psi)}(i)
=\prod_{i=2}^{k+\ell}c_{\vec{\rind},\vec{\sind}}^{\vec{\tind},(\varphi,\psi)}(i)
=\prod_{i=1}^{k+\ell-1}c_{\vec{\rind}\,',\vec{\sind}}^{\vec{\tind}\,',(\dt{\varphi},\st{\psi})}(i)
=
c_{\vec{\rind}\,',\vec{\sind}}^{\vec{\tind}\,',(\dt{\varphi},\st{\psi})}.
$$ This proves Eq.~(\mref{eq:coef21}). The proof of
Eq.~(\mref{eq:coef22}) is similar. The equality for $\tind_1\geq 2$
follows from Lemma \mref{lem:coefzero2}.(\mref{it:coefzero1}).

\smallskip

\noindent (\mref{it:coef3}) Suppose that $\rind_1=1$ and
$\sind_1\geq 2$.
\smallskip

\noindent
{\bf Case 1: $\tind_1=1$.} We consider the case of
$\varphi(1)=1$. By Lemma \mref{lem:coefind}.(\mref{it:coef1b}) we
have
$$c_{\vec{\rind},\vec{\sind}}^{\vec{\tind}, (\varphi,\psi)}(i+1)
=c_{\vec{\rind}\,',\vec{\sind}}^{\vec{\tind}\,',
(\dt{\varphi},\st{\psi})}(i).$$ Combining this with
$$c_{\vec{\rind}, \vec{\sind}}^{\vec{\tind}, (\varphi,
\psi)}(1)=\binc{\tind_1-1}{\rind_1-1}=\binc{0}{0}=1,$$ we obtain
$$c_{\vec{\rind},\vec{\sind}}^{\vec{\tind}, (\varphi, \psi)}
= \prod_{i=1}^{k+\ell} c_{\vec{\rind},\vec{\sind}}^{\vec{\tind},
(\varphi,\psi)}(i) =\prod_{i=1}^{k+\ell-1}
c_{\vec{\rind}\,',\vec{\sind}}^{\vec{\tind}\,',
(\dt{\varphi},\st{\psi})}(i)
=c_{\vec{\rind}\,',\vec{\sind}}^{\vec{\tind}\,',
(\dt{\varphi},\st{\psi})}.$$ This proves Eq.~(\mref{eq:coef31}). If
$\psi(1)=1$, then
$$ c_{\vec{\rind},\vec{\sind}}^{\vec{\tind},
(\varphi,\psi)}(1)=\binc{\tind_1-1}{\sind_1-1}=\binc{0}{\sind_1-1}=0
$$ since $\sind_1-1\geq 1$ and so
$ c_{\vec{\rind},\vec{\sind}}^{\vec{\tind},(\varphi,\psi)}=0,$ as needed.
\smallskip

\noindent
{\bf Case 2:  $\tind_1\geq 2$.}  We will consider the four subcases
when $\psi(1)=1$ and $\tind_1<\sind_1$, when $\psi(1)=1$ and
$\tind_1>\sind_1$, when $\psi(1)=1$ and $\tind_1=\sind_1$, and when
$\varphi(1)=1$.

If  $\psi(1)=1$ and $\tind_1<\sind_1$, then
$$
c_{\vec{\rind},\vec{\sind}}^{\vec{\tind},(\varphi,\psi)}=0=
c_{\vec{\rind},\vec{\sind}-\vec{e}_1}^{\vec{\tind}-\vec{e}_1,(\varphi,\psi)}$$
by Lemma \mref{lem:coefzero}. If $\psi(1)=1$ and $\tind_1>\sind_1$,
then by Lemma \mref{lem:coefzero2}.(\mref{it:coefzero1}) we also
have
$$
c_{\vec{\rind},\vec{\sind}}^{\vec{\tind},(\varphi,\psi)}=0=
c_{\vec{\rind},\vec{\sind}-\vec{e}_1}^{\vec{\tind}-\vec{e}_1,(\varphi,\psi)}.$$
So in these two subcases (\mref{eq:coef32}) holds.

Now if $\psi(1)=1$ and $\tind_1=\sind_1$, then
$$c_{\vec{\rind},\vec{\sind}}^{\vec{\tind},(\varphi,\psi)}(1)=\binc{\tind_1-1}{\sind_1-1}=
1= \binc{\tind_1-2}{\sind_1-2}
=c_{\vec{\rind},\vec{\sind}-\vec{e}_1}^{\vec{\tind}-\vec{e}_1,(\varphi,\psi)}(1).
$$ If $\varphi(1)=1$, then since $\rind_1=1$, we have
$$c^{\vec{\tind},
(\varphi,\psi)}_{\vec{\rind},
\vec{\sind}}(1)=\binc{\tind_1-1}{\rind_1-1} =\binc{\tind_1-1}{0}=1
=\binc{\tind_1-2}{0}=
\binc{\tind_1-2}{\rind_1-1}=c^{\vec{\tind}-\vec{e}_1,
(\varphi,\psi)}_{\vec{\rind}, \vec{\sind}-\vec{e}_1}(1).$$ In both
subcases, by Lemma \mref{lem:coefind}.(\mref{it:coef1a}) we always have
$$
c^{\vec{\tind},(\varphi,\psi)}_{\vec{\rind},\vec{\sind}}(i)=c^{\vec{\tind}-\vec{e}_1,
(\varphi,\psi)}_{\vec{\rind}, \vec{\sind}-\vec{e}_1}(i).
$$
for $i\geq 2$. Therefore,
$$
c^{\vec{\tind},(\varphi,\psi)}_{\vec{\rind},\vec{\sind}}
=\prod_{i=1}^{k+\ell}c^{\vec{\tind},(\varphi,\psi)}_{\vec{\rind},\vec{\sind}}(i)
=\prod_{i=1}^{k+\ell}c^{\vec{\tind}-\vec{e}_1,
(\varphi,\psi)}_{\vec{\rind},\vec{\sind}-\vec{e}_1}(i)
=c^{\vec{\tind}-\vec{e}_1,(\varphi,\psi)}_{\vec{\rind},\vec{\sind}-\vec{e}_1}
.$$ This proves (\mref{eq:coef32}).

The proof for the instance of $\rind_1\geq 2$ and $\sind_1=1$ is
similar.
\end{proof}

\section{Proof of the main theorems}
\mlabel{sec:proof}
We first show that, under the condition that $\sg$ is an abelian group,  Theorem~\mref{thm:qshsh} and Theorem~\mref{thm:qshshe} are equivalent. Then we only need to prove
Theorem~\mref{thm:qshsh}. This is done in Section~\mref{ss:proof}.

\subsection{The equivalence between Theorem~\mref{thm:qshsh} and Theorem~\mref{thm:qshshe}}
\mlabel{ss:equiv}
We start with a lemma.
\begin{lemma} Let $\sg$ be an abelian group. With the notations in Eq.~(\mref{eq:la1}), (\mref{eq:mulind}) and (\mref{eq:pssha}), we have
\begin{equation} \theta(\vec{a}\ssha_{(\varphi,\psi)}\vec{b}) =
\theta(\vec{a})\pssha_{(\varphi,\psi)}\theta(\vec{b}).
\mlabel{eq:psheqsh}
\end{equation}
\mlabel{lem:psheqsh}
\end{lemma}
\begin{proof}
Let $\vec{w}=\theta(\vec{a})$ and
$\vec{z}=\theta(\vec{b})$. Then by Eq.~(\mref{eq:la1}), we have $w_j=1/a_1$ when $j=1$ and $w_j=a_{j-1}/a_j$ when $j\geq 2$.
Similarly, $z_j=1/b_1$ when $j=1$ and $z_j=b_{j-1}/b_j$ when $j\geq 2$.

Recall Eq.~(\mref{eq:mulind}):
$$(\vec{a}\ssha_{(\varphi,\psi)}\vec{b})_i
=\left\{
\begin{array}{ll} a_j & \text{ if } i=\varphi(j),
\\
b_j & \text{ if } i=\psi(j).
\end{array}
\right.$$
When $i=1$, we have
$$\theta(\vec{a}\ssha_{(\varphi,\psi)}\vec{b})_1=
(\vec{a}\ssha_{(\varphi,\psi)}\vec{b})_1^{-1} =
\left\{\begin{array}{ll} a_1^{-1}=w_1 & \text{ if }1=\varphi(1) \\
b_1^{-1}=z_1 & \text{ if }1=\psi(1)
\end{array}\right. = (\vec{w}\pssha_{(\varphi,\psi)} \vec{z})_1.
$$ Next let $i\geq 2$. Assume that $i\in\im(\varphi)$, say $i=\varphi(j)$ for some $j\in [k]$.
If $i-1\in \im(\varphi)$, then $j\geq 2$ and $i-1=\varphi(j-1)$.
Thus
$$\theta(\vec{a}\ssha_{(\varphi,\psi)}\vec{b})_i
=\frac{(\vec{a}\ssha_{(\varphi,\psi)}\vec{b})_{i-1}}{
(\vec{a}\ssha_{(\varphi,\psi)}\vec{b})_i} = \frac{a_{j-1}}{a_j} =w_j.
$$ If $i-1\in \im(\psi)$, then $i-1=\psi(i-j)$. Thus
$$\theta(\vec{a}\ssha_{(\varphi,\psi)}\vec{b})_i
=\frac{(\vec{a}\ssha_{(\varphi,\psi)}\vec{b})_{i-1}}{
(\vec{a}\ssha_{(\varphi,\psi)}\vec{b})_i} = \frac{b_{i-j}}{a_j}=\frac{(z_1\cdots z_{i-j})^{-1}}{(w_1\cdots w_j)^{-1}}=\frac{w_1\cdots
w_j}{z_1\cdots z_{i-j}}.
$$ Hence by Eq.~(\mref{eq:pssha}),
$$\theta(\vec{a}\ssha_{(\varphi,\psi)}\vec{b})_i=(\vec{w}\pssha_{(\varphi,\psi)}\vec{z})_i$$
when $i\in \im(\varphi)$. A similar argument shows that the above
equality also holds when $i\in \im(\psi)$. This proves
(\mref{eq:psheqsh}).
\end{proof}

\begin{prop}
When $\sg$ is an abelian group,
Theorem~\mref{thm:qshshe} is equivalent to Theorem~\mref{thm:qshsh}.
\mlabel{pp:equiv}
\end{prop}
\begin{proof}
{}From the definitions of $\theta$, $\qsshae$ and $\qsshab$, we see that $\theta$ is an algebra isomorphism from
$\calh\bsh(\zsg{\sg})=(\calh(\zsg{\sg}),\qsshab)$ to
$\calh\esh(\zsg{\sg})=(\calh(\zsg{\sg}),\qsshae)$.
So for any $\spair{\vec{r}}{\vec{a}}, \spair{\vec{s}}{\vec{b}}\in \calh\bsh(\zsg{\sg})$,
$$
\begin{aligned}
&\spair{\vec{r}}{\vec{a}} \qsshab \spair{\vec{s}}{\vec{b}}
= \sum_{\substack{ (\varphi,\psi)\in \indI_{k,\ell}\\
\vec{\tind}\in \ZZ_{\geq 1}^{k+\ell},
|\vec{\tind}|=|\vec{\rind}|+|\vec{\sind}| }}
c_{\vec{\rind},\vec{\sind}}^{\vec{\tind},
  (\varphi,\psi)}
  \spair{\vec{\tind}}{\vec{a}\ssha_{(\varphi,\psi)}\vec{b}}
\\
\Leftrightarrow \hspace{.2cm}
& \theta(\spair{\vec{r}}{\vec{a}} \qsshab \spair{\vec{s}}{\vec{b}})
=
\theta( \sum_{\substack{ (\varphi,\psi)\in \indI_{k,\ell}\\
\vec{\tind}\in \ZZ_{\geq 1}^{k+\ell},
|\vec{\tind}|=|\vec{\rind}|+|\vec{\sind}| }}
c_{\vec{\rind},\vec{\sind}}^{\vec{\tind},
  (\varphi,\psi)}
  \spair{\vec{\tind}}{\vec{a}\ssha_{(\varphi,\psi)}\vec{b}})\\
\Leftrightarrow \hspace{.2cm}
& \theta(\spair{\vec{r}}{\vec{a}}) \qsshae \theta (\spair{\vec{s}}{\vec{b}})
=
 \sum_{\substack{ (\varphi,\psi)\in \indI_{k,\ell}\\
\vec{\tind}\in \ZZ_{\geq 1}^{k+\ell},
|\vec{\tind}|=|\vec{\rind}|+|\vec{\sind}| }}
c_{\vec{\rind},\vec{\sind}}^{\vec{\tind},
  (\varphi,\psi)}
  \theta(\spair{\vec{\tind}}{\vec{a}\ssha_{(\varphi,\psi)}\vec{b}})
\\
\Leftrightarrow \hspace{.2cm}
& \spair{\vec{r}}{\theta(\vec{a})} \qsshae \spair{\vec{s}}{\theta (\vec{b})}
=
 \sum_{\substack{ (\varphi,\psi)\in \indI_{k,\ell}\\
\vec{\tind}\in \ZZ_{\geq 1}^{k+\ell},
|\vec{\tind}|=|\vec{\rind}|+|\vec{\sind}| }}
c_{\vec{\rind},\vec{\sind}}^{\vec{\tind},
  (\varphi,\psi)}
  \spair{\vec{\tind}}{\theta(\vec{a}\ssha_{(\varphi,\psi)}\vec{b})}
\\
\Leftrightarrow \hspace{.2cm}
& \spair{\vec{r}}{\theta(\vec{a})} \qsshae \spair{\vec{s}}{\theta (\vec{b})}
=
 \sum_{\substack{ (\varphi,\psi)\in \indI_{k,\ell}\\
\vec{\tind}\in \ZZ_{\geq 1}^{k+\ell},
|\vec{\tind}|=|\vec{\rind}|+|\vec{\sind}| }}
c_{\vec{\rind},\vec{\sind}}^{\vec{\tind},
  (\varphi,\psi)}
  \spair{\vec{\tind}}{\theta(\vec{a})\pssha_{(\varphi,\psi)}\theta(\vec{b})}
    \qquad (\text{by Eq.~(\mref{eq:psheqsh})}).
  \end{aligned}
$$
Then the proposition follows from the bijectivity of $\theta$.
\end{proof}

\subsection{Proof of Theorem~\mref{thm:qshsh}}
\mlabel{ss:proof}

In this section we prove Theorem~\mref{thm:qshsh}. We first describe recursive relations of
$\qsshab$ that we will use later in the proof.

Let $\calh\bsh^+(\zsg{\sg})$ be the subring of
$\calh\bsh(\zsg{\sg})$ generated by $\spair{\vec{\sind}}{\vec{b}}$
with $\vec{\sind}\in \ZZ_{\geq 1}^k, \vec{b}\in\sg^k, k\geq 1$. Then
$$ \calh\bsh(\zsg{\sg})= \ZZ\oplus \calh\bsh^+(\zsg{\sg}). $$
Define the following operators
$$
\begin{aligned}
P:& \calh\bsh^+(\zsg{\sg})\rightarrow \calh\bsh(\zsg{\sg}), \quad
P(\spair{\sind_1,\sind_2,\cdots ,\sind_k}{b_1,b_2,\cdots,b_k})
= \spair{\sind_1+1, \sind_2, \cdots, \sind_k}{\;\;\;\;\;b_1,b_2,\cdots,b_k}, \\
Q_b:& \calh\bsh(\zsg{\sg})\rightarrow \calh\bsh(\zsg{\sg}), \quad
Q_b(\spair{\sind_1,\cdots, \sind_k}{b_1,\cdots, b_k})= \spair{ 1,
 \sind_1, \cdots, \sind_k}{b,b_1, \cdots, b_k }, \quad
Q_b(1)=\spair{1}{b}.
\end{aligned}
$$

\begin{prop} The multiplication $\qsshab$ on
$\calh\bsh(\zsg{\sg})$ defined in Eq.~(\mref{eq:shtrans}) is the
unique one that satisfies the Rota-Baxter type relations~\mcite{G-K1}:
$$
\begin{aligned}
&P(\xi_1)\qsshab P(\xi_2) = P\big(\xi_1 \qsshab P(\xi_2)\big)+
P\big(P(\xi_1)\qsshab \xi_2\big),  \  \xi_1,\xi_2\in\calh\bsh^+(\zsg{\sg}), \\
&Q_{a}(\xi_1)\qsshab Q_{b}(\xi_2) = Q_{a}\big(\xi_1 \qsshab
Q_{b}(\xi_2)\big)
+ Q_{b}\big( Q_{a}(\xi_1)\qsshab \xi_2\big),  \  \xi_1,\xi_2\in\calh\bsh(\zsg{\sg}), \\
&P(\xi_1)\qsshab Q_{b}(\xi_2) = Q_{b}\big(P(\xi_1)\qsshab
\xi_2\big)+ P\big(\xi_1\qsshab Q_{b}(\xi_2)\big), \
\xi_1\in\calh\bsh^+(\zsg{\sg}),
\xi_2 \in\calh\bsh(\zsg{\sg}), \\
&Q_{b}(\xi_1)\qsshab P(\xi_2) = Q_{b}\big(\xi_1\qsshab
P(\xi_2)\big)+ P\big(Q_{b}(\xi_1)\qsshab \xi_2\big), \
\xi_1\in\calh\bsh(\zsg{\sg}), \xi_2 \in\calh\bsh^+(\zsg{\sg}).
\end{aligned}
$$
with the initial condition that $1 \qsshab \xi=\xi\qsshab 1 = \xi$
for $\xi\in \calh\bsh(\zsg{\sg})$. \mlabel{pp:rec}
\end{prop}
\begin{proof}
Let $\calh\shf\lone^+(\ssg{\sg})$ be the subring of
$\calh\shf\lone(\ssg{\sg})$ generated by words of the form
$u x_b$ with $b\in \sg$. %  $(w\in\sg)$.
Then
$$ \calh\shf\lone(\ssg{\sg})= \ZZ\oplus \calh\shf\lone^+(\ssg{\sg}). $$
Define operators
$$
\begin{aligned}
I_0: & \calh\shf\lone^+(\ssg{\sg})\to \calh\shf\lone(\ssg{\sg}), \quad I_0(u)= x_0 u,\\
I_b: & \calh\shf\lone(\ssg{\sg}) \to \calh\shf\lone(\ssg{\sg}),
\quad I_b(u)=\left \{\begin{array}{ll} x_b u, & u \neq 1, \\
x_b, & u=1,
\end{array}\right .
\end{aligned}
$$
 for $b\in\sg$.
Then the well-known recursive formula of the shuffle product
$$ (a_1 \fraka) \ssha (b_1 \frakb) = a_1  (\fraka \ssha (b_1  \frakb))+b_1  ((a_1 \fraka)\ssha \frakb), a_1,b_1\in \ssg{\sg}, \fraka,\frakb \in M (\ssg{\sg})$$
can be rewritten as the
following relations of $I_0$ and $I_a$, $I_b,$ $a,b\in\sg$,
\begin{equation}
\begin{aligned}
I_0(u) \ssha I_0(v) &= I_0\big( u\ssha I_0(v))+ I_0(I_0(u)\ssha v\big),
\quad u,v\in \calh\shf\lone^+(\ssg{\sg}),
\\
I_a(u) \ssha I_b(v) &= I_a\big( u\ssha I_b(v))+ I_b(I_a(u)\ssha
v\big), \quad u,v\in \calh\shf\lone(\ssg{\sg}),
\\
I_0(u) \ssha I_b(v) &= I_0\big( u\ssha I_b(v)\big)+
I_b\big(I_0(u)\ssha v\big), \quad  u\in \calh\shf\lone^+(\ssg{\sg}),
v\in \calh\shf\lone(\ssg{\sg}),
\\
I_b(u) \ssha I_0(v) &= I_b\big( u\ssha I_0(v)\big)+
I_0\big(I_b(u)\ssha v\big), \quad  u\in
\calh\shf\lone(\ssg{\sg}),v\in \calh\shf\lone^+(\ssg{\sg}).
\end{aligned}
\mlabel{eq:iform}
\end{equation}
Under the bijection $\rho: \calh\shf\lone(\ssg{\sg})\to
\calh\bsh(\zsg{\sg})$ in Eq.~(\mref{eq:shqsh1}), $I_0$ and $I_b,$
$b\in\sg$, are sent to $P$ and $Q_{b},$ $b\in\sg$, respectively.
Further the relations in Eq.~(\mref{eq:iform}) for $I_0$ and $I_b,$
$b\in\sg$, take the form in Proposition~\mref{pp:rec}. Finally,
since $\ssha$ is the unique multiplication on
$\calh\shf\lone(\ssg{\sg})$ characterized by its recursive relation
Eq.~(\mref{eq:iform}) and the initial condition $1 \ssha u= u\ssha
1=u$, $\qsshab$ is also unique as characterized.
\end{proof}

For $\vec{b}\in\sg^k$, recall the following notation from Definition~\mref{de:notn}:
$$\vec{b}\,'=(b'_1,\cdots, b'_{k-1})
:=(b_2,\cdots, b_k)$$
 with the convention that
$\vec{b}\,' = \bfe$ when $k=1$.
In the proof for Theorem \mref{thm:qshsh} we also need the following
lemma.

\begin{lemma} Let $\vec{\tind}\in \ZZ_{\geq 1}^{k+\ell-1}$,
$\vec{a}\in\sg^k$ and $\vec{b}\in\sg^\ell$.
\begin{enumerate}
\item
For any $(\varphi,\psi)\in\indI_{k-1,\ell}$ we have
\begin{equation}
Q_{a_1}(\spair{\vec{\tind}}{ \vec{a}\,'\ssha_{(\varphi,\psi)}\vec{b}})
= \spair{(1,\vec{\tind})} {\vec{a}\ssha_{(\varphi^{\&},
\psi^{*})}\vec{b}} \mlabel{eq:pregp1}
\end{equation} with the notations in Eq.~(\mref{eq:mulind}) and Definition~\mref{de:notn}.
\mlabel{it:pregp1}
\item
For any $(\varphi,\psi)\in\indI_{k,\ell-1}$ we have
\begin{equation}
Q_{b_1}(\spair{\vec{\tind}}{ \vec{a}\ssha_{(\varphi,\psi)}\vec{b}\,'})
= \spair{(1,\vec{\tind})}{\vec{a}\ssha_{(\varphi^{*},
\psi^{\&})}\vec{b}}. \mlabel{eq:pregp2}
\end{equation}
\mlabel{it:pregp2}
\end{enumerate}
\mlabel{lem:pre-gp}
\end{lemma}
\begin{proof} (\mref{it:pregp1}) Let $\vec{\varpi}=(\varpi_1,\cdots,\varpi_{k+\ell-1}):
=\vec{a}\,'\ssha_{(\varphi,\psi)}\vec{b}$ and
$\vec{\tau}=(\tau_1,\cdots,\tau_{k+\ell}):=
\vec{a}\ssha_{(\varphi^{\&},\psi^*)}\vec{b}$. By the definition of
$Q_{a_1}$, we only need to prove that
$$\tau_i =\left\{
                  \begin{array}{ll}
                  a_1 & \text{ if } i=1, \\
                  \varpi_{i-1} & \text{ if } i\geq 2.
                  \end{array}
             \right. $$
Since $\varphi^{\&}(1)=1$, we have $\tau_1=a_1$. Now let $i\geq 2$.
We have $i\in \im(\varphi^{\&})$ or $i\in \im(\psi^*)$. If $i\in
\im(\varphi^{\&})$, say $i=\varphi^{\&}(j)$, then
$i-1=\varphi(j-1)$. Thus we have $\tau_i= a_j$ and
$\varpi_{i-1}=a'_{j-1}=a_j$. This shows that $\tau_i=\varpi_{i-1}$.
If $i\in \im(\psi^*)$, say $i=\psi^*(j)$, then $i-1=\psi(j)$. Thus
$\tau_i=b_j$ and $\varpi_{i-1}=b_j$ again showing
$\tau_i=\varpi_{i-1}$.

(\mref{it:pregp2}). The proof is similar to that for
Item.~(\mref{it:pregp1}).
\end{proof}

\medskip

\noindent{\it Proof of Theorem \mref{thm:qshsh}}. We prove the
extended form of (\mref{eqn:maincoef}) where one of $k$ and $\ell$,
but not both, might be zero. We prove this by induction on
$|\vec{\rind}|+|\vec{\sind}|\geq 1$. If
$|\vec{\rind}|+|\vec{\sind}|= 1$, then exactly one of $k$ and $\ell$
is zero. So exactly one of $\spair{\vec{\rind}}{\vec{a}}$ and
$\spair{\vec{\sind}}{\vec{b}}$ is the identity $1$. Then by
(\mref{eq:memp}) and (\mref{eq:nemp}), there is nothing to prove.
For any given integer $n\geq 2$, assume that the assertion holds for
every pair $(\vec{\rind},\vec{\sind})$ with
$|\vec{\rind}|+|\vec{\sind}|<n$. Now consider $\vec{\rind}$ and
$\vec{\sind}$ with $|\vec{\rind}|+|\vec{\sind}|=n$.  If one of $k$
or $\ell$ is $0$, then again by (\mref{eq:memp}) and
(\mref{eq:nemp}) there is nothing to prove. So we may assume that
$k,\ell\geq 1$. There are four cases to consider.
\smallskip

\noindent {\bf Case 1. $\rind_1\geq 2$ and $\sind_1\geq 2$.}  Then
by Proposition~\mref{pp:rec} and the induction hypothesis, we have
{\allowdisplaybreaks
\begin{eqnarray*}
\spair{\vec{\rind}}{\vec{a}}\qsshab \spair{\vec{\sind}}{\vec{b}} &=&
P(\spair{\vec{\rind}-\vec{e}_1}{\vec{a}})\qsshab
P(\spair{\vec{\sind}-\vec{e}_1}{\vec{b}})
    \\
&=& P(\spair{\vec{\rind}-\vec{e}_1}{\vec{a}}\qsshab
    \spair{\vec{\sind}}{\vec{b}}+ \spair{\vec{\rind}}{\vec{a}}\qsshab \spair{\vec{\sind}-\vec{e}_1}{\vec{b}})
    \\
&=&
P\big(\sum_{(\varphi,\psi)\in\indI_{k,\ell}} \sum_{\substack{\vec{\tind}\:\in
    \: \ZZ_{\geq 1}^{k+\ell}\\
    |\vec{\tind}\:|=|\vec{\rind}|+|\vec{\sind}|-1}}
    (c^{\vec{\tind},(\varphi,\psi)}_{\vec{\rind}-\vec{e}_1,\vec{\sind}}
    +c^{\vec{\tind},(\varphi,\psi)}_{\vec{\rind},\vec{\sind}-\vec{e}_1})\spair{\vec{\tind}}{
    \vec{a}\ssha_{(\varphi,\psi)}\vec{b}}\,\big)
    \\
&=& \sum_{(\varphi,\psi)\in\indI_{k,\ell}} \sum_{\substack{\vec{\tind}\:\in \:
     \ZZ_{\geq 1}^{k+\ell}\\
    |\vec{\tind}\:|=|\vec{\rind}|+|\vec{\sind}|-1}}
    (c^{\vec{\tind},(\varphi,\psi)}_{\vec{\rind}-\vec{e}_1,\vec{\sind}}
    +c^{\vec{\tind},(\varphi,\psi)}_{\vec{\rind},\vec{\sind}-\vec{e}_1})
    \spair{\vec{\tind}+\vec{e}_1}{\vec{a}\ssha_{(\varphi,\psi)}\vec{b}}
    \\
&=& \sum_{(\varphi,\psi)\in\indI_{k,\ell}} \sum_{\substack{\vec{\tind}\:\in \:
    \ZZ_{\geq 1}^{k+\ell}\\
    |\vec{\tind}\:|=|\vec{\rind}|+|\vec{\sind}|, \tind_1\geq 2}}
    (c^{\vec{\tind}-\vec{e}_1,(\varphi,\psi)}_{\vec{\rind}-\vec{e}_1,\vec{\sind}}
    +c^{\vec{\tind}-\vec{e}_1,(\varphi,\psi)}_{\vec{\rind},\vec{\sind}-\vec{e}_1})
    \spair{\vec{\tind}}{
    \vec{a}\ssha_{(\varphi,\psi)}\vec{b}}
    \\
&=& \sum_{(\varphi,\psi)\in\indI_{k,\ell}} \sum_{\substack{\vec{\tind}\:\in \:
    \ZZ_{\geq 1}^{k+\ell}\\
    |\vec{\tind}\:|=|\vec{\rind}|+|\vec{\sind}|, \tind_1\geq 2}}
    c^{\vec{\tind},(\varphi,\psi)}_{\vec{\rind},\vec{\sind}}\spair{\vec{\tind}}{
    \vec{a}\ssha_{(\varphi,\psi)}\vec{b}}
    \qquad \text{(by Eq.~(\mref{eq:coef1}))}
    \\
&=& \sum_{(\varphi,\psi)\in\indI_{k,\ell}} \sum_{\substack{\vec{\tind}\:\in \:
    \ZZ_{\geq 1}^{k+\ell}\\
    |\vec{\tind}\:|=|\vec{\rind}|+|\vec{\sind}|}}
    c^{\vec{\tind},(\varphi,\psi)}_{\vec{\rind},\vec{\sind}}\spair{\vec{\tind}}{
    \vec{a}\ssha_{(\varphi,\psi)}\vec{b}}
    \qquad \text{(by Eq.~(\mref{eq:coefzero1}))}.
\end{eqnarray*}
}

\noindent {\bf Case 2.  $\rind_1=\sind_1=1$.} We will use the
notations $\vec{\rind}\,',\vec{\sind}\,', \vec{a}\,'$ and $\vec{b}\,'$ in
Definitions~\mref{de:notn}. Then
{\allowdisplaybreaks
\begin{eqnarray*}
\lefteqn{\spair{\vec{\rind}}{\vec{a}}\qsshab
\spair{\vec{\sind}}{\vec{b}} =
     Q_{a_1}(\spair{\vec{\rind}\,'}{\vec{a}\,'})\qsshab Q_{b_1}(\spair{\vec{\sind}\,'}{\vec{b}\,'})}
     \\
&=&  Q_{a_1}(\spair{\vec{\rind}\,'}{\vec{a}\,'}\qsshab
     \spair{\vec{\sind}}{\vec{b}})+Q_{b_1}
     (\spair{\vec{\rind}}{\vec{a}}\qsshab \spair{\vec{\sind}\,'}{\vec{b}\,'})
     \\
&=&
     Q_{a_1}\big(\sum_{(\varphi,\psi)\in\indI_{k-1,\ell}} \sum_{\substack{\vec{\tind}\in
      \:\ZZ_{\geq 1}^{k+\ell-1}\\ |\vec{\tind}|=|\vec{\rind}|+|\vec{\sind}|-1}}
      c_{\vec{\rind}\,',\vec{\sind}}^{\vec{\tind},(\varphi,\psi)}
      \spair{\vec{\tind}}{\vec{a}\,'\ssha_{(\varphi,\psi)}\vec{b}}\big)
          \\
          && \quad + Q_{b_1}\big(\sum_{(\varphi,\psi)\in\indI_{k,\ell-1}} \sum_{\substack{\vec{\tind}\in \:\ZZ_{\geq
      1}^{k+\ell-1}\\ |\vec{\tind}|=|\vec{\rind}|+|\vec{\sind}|-1}}
      c_{\vec{\rind},\vec{\sind}\,'}^{\vec{\tind},
      (\varphi,\psi)}\spair{\vec{\tind}}{
      \vec{a}\ssha_{(\varphi,\psi)}\vec{b}\,'}\,\big)
      \\
&=&  \sum_{(\varphi,\psi)\in\indI_{k-1,\ell}} \sum_{\substack{\vec{\tind}\in \:
      \ZZ_{\geq 1}^{k+\ell-1}\\ |\vec{\tind}|=|\vec{\rind}|+|\vec{\sind}|-1}}
      c_{\vec{\rind}\,',\vec{\sind}}^{\vec{\tind},(\varphi,\psi)}
      \spair{(1,\vec{\tind})}{
      \vec{a}\ssha_{(\varphi^{\&},\psi^*)}\vec{b}}
      \\
      && \quad + \sum_{(\varphi,\psi)\in\indI_{k,\ell-1}} \sum_{\substack{\vec{\tind}\in \: \ZZ_{\geq 1}^{k+\ell-1}\\
      |\vec{\tind}|=|\vec{\rind}|+|\vec{\sind}|-1}}
      c_{\vec{\rind},\vec{\sind}\,'}^{\vec{\tind},(\varphi,\psi)}
      \spair{(1,\vec{\tind})}{
      \vec{a}\ssha_{(\varphi^{*},\psi^{\&})}\vec{b}}
      \quad (\text{ by Eq.}
      (\mref{eq:pregp1}) \text{ and }(\mref{eq:pregp2}))
      \\
&=&
      \sum_{\substack{(\varphi,\psi)\in\indI_{k,\ell}\\ \varphi(1)=1}} \sum_{\substack{\vec{\tind}\in
      \: \ZZ_{\geq 1}^{k+\ell-1}\\ |\vec{\tind}|=|\vec{\rind}|+|\vec{\sind}|-1}}
      c_{\vec{\rind}\,',\vec{\sind}}^{\vec{\tind},(\dt{\varphi},\st{\psi})}
      \spair{(1,\vec{\tind})}{
      \vec{a}\ssha_{(\varphi,\psi)}\vec{b}}
      \\
      && \quad + \sum_{\substack{(\varphi,\psi)\in\indI_{k,\ell}\\ \psi(1)=1}}\sum_{\substack{\vec{\tind}\in \: \ZZ_{\geq 1}^{k+\ell-1}\\
      |\vec{\tind}|=|\vec{\rind}|+|\vec{\sind}|-1}}
      c_{\vec{\rind},\vec{\sind}\,'}^{\vec{\tind},(\st{\varphi},\dt{\psi})}
      \spair{(1,\vec{\tind})}{
      \vec{a}\ssha_{(\varphi,\psi)}\vec{b}} \quad \text{(by Lemma~\mref{lem:bij})}
      \\
&=&
      \sum_{\substack{(\varphi,\psi)\in\indI_{k,\ell}\\ \varphi(1)=1}} \sum_{\substack{\vec{\tind}\in
      \: \ZZ_{\geq 1}^{k+\ell-1}\\ |\vec{\tind}|=|\vec{\rind}|+|\vec{\sind}|-1}}
      c_{\vec{\rind},\vec{\sind}}^{(1,\vec{\tind}),(\varphi,\psi)}
      \spair{(1,\vec{\tind})}{
      \vec{a}\ssha_{(\varphi,\psi)}\vec{b}}
      \\
      && \quad + \sum_{\substack{(\varphi,\psi)\in\indI_{k,\ell}\\ \psi(1)=1}} \sum_{\substack{\vec{\tind}\in \: \ZZ_{\geq 1}^{k+\ell-1}\\
      |\vec{\tind}|=|\vec{\rind}|+|\vec{\sind}|-1}}
      c_{\vec{\rind},\vec{\sind}}^{(1,\vec{\tind}),(\varphi,\psi)}
      \spair{(1,\vec{\tind})}{
      \vec{a}\ssha_{(\varphi,\psi)}\vec{b}} \quad \text{(by Eq.~(\mref{eq:coef21}) and (\mref{eq:coef22}))}
      \\
&=&   \sum_{(\varphi,\psi)\in\indI_{k,\ell}} \sum_{\substack{\vec{\tind}\in \:
      \ZZ_{\geq 1}^{k+\ell-1}\\
      |\vec{\tind}|=|\vec{\rind}|+|\vec{\sind}|-1}}
      c_{\vec{\rind},\vec{\sind}}^{(1,\vec{\tind}),(\varphi,\psi)}
      \spair{(1,\vec{\tind})}{
      \vec{a}\ssha_{(\varphi,\psi)}\vec{b}}
      \\
&=&
      \sum_{(\varphi,\psi)\in\indI_{k,\ell}} \sum_{\substack{\vec{\tind}\in\ZZ_{\geq
      1}^{k+\ell}\\ |\vec{\tind}|=|\vec{\rind}|+|\vec{\sind}|,\tind_1=1}}
      c_{\vec{\rind},\vec{\sind}}^{\vec{\tind},(\varphi,\psi)}
      \spair{\vec{\tind}}{
      \vec{a}\ssha_{(\varphi,\psi)}\vec{b}}
      \\
&=&
      \sum_{(\varphi,\psi)\in\indI_{k,\ell}} \sum_{\substack{\vec{\tind}\in
      \:\ZZ_{\geq 1}^{k+\ell}\\ |\vec{\tind}\:|=|\vec{\rind}|+|\vec{\sind}|}}
           c_{\vec{\rind},\vec{\sind}}^{\vec{\tind},(\varphi,\psi)}
           \spair{\vec{\tind}}{\vec{a}\ssha_{(\varphi,\psi)}\vec{b}}
       \qquad \text{ (by Eq.~(\mref{eq:coefzero2}))}.
\end{eqnarray*}
}

\noindent {\bf Case 3. $\rind_1=1$ and $\sind_1\geq 2$.} With the
notations in Definitions~\mref{de:notn}, we
write $\vec{\rind}=(1,\vec{\rind}\,')$. Let
$\vec{a}\,'=(a_2,\cdots, a_r)$. Then
{\allowdisplaybreaks
\begin{eqnarray*}
\lefteqn{ \spair{\vec{\rind}}{\vec{a}}\qsshab
\spair{\vec{\sind}}{\vec{b}} =
Q_{a_1}(\spair{\vec{\rind}\,'}{\vec{a}\,'})\qsshab
P(\spair{\vec{\sind}-\vec{e}_1}{\vec{b}})}
\\
 &=& Q_{a_1}(\spair{\vec{\rind}\,'}{\vec{a}\,'}\qsshab \spair{\vec{\sind}}{\vec{b}})+
 P(\spair{\vec{\rind}}{\vec{a}}\qsshab
     \spair{\vec{\sind}-\vec{e}_1}{\vec{b}})
 \\
 &=& Q_{a_1} (\sum_{(\varphi,\psi)\in\indI_{k-1,\ell}} \sum_{\substack{\vec{\tind}\in \:\ZZ_{\geq 1}^{k+\ell-1}\\
       |\vec{\tind}|=|\vec{\rind}|+|\vec{\sind}|-1}}
       c_{\vec{\rind}\,',\vec{\sind}}^{\vec{\tind},(\varphi,\psi)}
       \spair{\vec{\tind}}{
       \vec{a}\,'\ssha_{(\varphi,\psi)}\vec{b}})
       \\
 && \quad + P(\sum_{(\varphi,\psi)\in\indI_{k,\ell}} \sum_{\substack{\vec{\tind}\in \:\ZZ_{\geq 1}^{k+\ell}\\
       |\vec{\tind}|=|\vec{\rind}|+|\vec{\sind}|-1}}
       c_{\vec{\rind},\vec{\sind}-\vec{e}_1}^{\vec{\tind},(\varphi,\psi)}
       \spair{\vec{\tind}}{\vec{a}\ssha_{(\varphi,\psi)}\vec{b}})
  \\
 &=& \sum_{(\varphi,\psi)\in\indI_{k-1,\ell}} \sum_{\substack{\vec{\tind}\in \:\ZZ_{\geq 1}^{k+\ell-1}\\
       |\vec{\tind}|=|\vec{\rind}|+|\vec{\sind}|-1}}
       c_{\vec{\rind}\,',\vec{\sind}}^{\vec{\tind},(\varphi,\psi)}
       \spair{(1,\vec{\tind})}{\vec{a}\ssha_{(\varphi^{\&},\psi^*)}\vec{b}}
       \\
  && \quad   + \sum_{(\varphi,\psi)\in\indI_{k,\ell}} \sum_{\substack{\vec{\tind}\in \:\ZZ_{\geq 1}^{k+\ell}\\
       |\vec{\tind}|=|\vec{\rind}|+|\vec{\sind}|-1}}
       c_{\vec{\rind},\vec{\sind}-\vec{e}_1}^{\vec{\tind}, (\varphi,\psi)}
       \spair{\vec{\tind}+\vec{e}_1}{\vec{a}\ssha_{(\varphi,\psi)}\vec{b}}
       \quad (\text{by Eq.}~ (\mref{eq:pregp1}))
 \\
 &=& \sum_{\substack{(\varphi,\psi)\in\indI_{k,\ell}\\\varphi(1)=1}} \sum_{\substack{\vec{\tind}\in \:\ZZ_{\geq 1}^{k+\ell-1}\\
       |\vec{\tind}|=|\vec{\rind}|+|\vec{\sind}|-1}}
       c_{\vec{\rind}\,',\vec{\sind}}^{\vec{\tind},(\dt{\varphi},\st{\psi})}
       \spair{(1,\vec{\tind})}{\vec{a}\ssha_{(\varphi,\psi)}\vec{b}}
       \\
  && \quad   + \sum_{(\varphi,\psi)\in\indI_{k,\ell}} \sum_{\substack{\vec{\tind}\in \:\ZZ_{\geq 1}^{k+\ell}\\
       |\vec{\tind}|=|\vec{\rind}|+|\vec{\sind}|-1}}
       c_{\vec{\rind},\vec{\sind}\,'}^{\vec{\tind},(\varphi,\psi)}
       \spair{\vec{\tind}+\vec{e}_1}{\vec{a}\ssha_{(\varphi,\psi)}\vec{b}}
       \quad \text{(by Lemma~\mref{lem:bij})}
 \\
 &=& \sum_{\substack{(\varphi,\psi)\in\indI_{k,\ell}\\ \varphi(1)=1}}\sum_{\substack{\vec{\tind}\in \:\ZZ_{\geq 1}^{k+\ell-1}\\
       |\vec{\tind}|=|\vec{\rind}|+|\vec{\sind}|-1}}
       c_{\vec{\rind},\vec{\sind}}^{(1,\vec{\tind}),(\varphi,\psi)}
       \spair{(1,\vec{\tind})}{\vec{a}\ssha_{(\varphi,\psi)}\vec{b}}
       \\
  && \quad   + \sum_{(\varphi,\psi)\in\indI_{k,\ell}} \sum_{\substack{\vec{\tind}\in \:\ZZ_{\geq 1}^{k+\ell}\\
       |\vec{\tind}|=|\vec{\rind}|+|\vec{\sind}|-1}}
       c_{\vec{\rind},\vec{\sind}}^{\vec{\tind}+\vec{e}_1,(\varphi,\psi)}
       \spair{\vec{\tind}+\vec{e}_1}{\vec{a}\ssha_{(\varphi,\psi)}\vec{b}}
       \quad (\text{by  Eq. ~ (\mref{eq:coef31}) and (\mref{eq:coef32})} )
  \\
  &=& \sum_{(\varphi,\psi)\in\indI_{k,\ell}} \sum_{\substack{\vec{\tind}\in \:\ZZ_{\geq 1}^{k+\ell-1}\\
       |\vec{\tind}|=|\vec{\rind}|+|\vec{\sind}|-1}}
       c_{\vec{\rind},\vec{\sind}}^{(1,\vec{\tind}),(\varphi,\psi)}
       \spair{(1,\vec{\tind})}{\vec{a}\ssha_{(\varphi,\psi)}\vec{b}}
       \\
  && \quad   + \sum_{(\varphi,\psi)\in\indI_{k,\ell}} \sum_{\substack{\vec{\tind}\in \:\ZZ_{\geq 1}^{k+\ell}\\
       |\vec{\tind}|=|\vec{\rind}|+|\vec{\sind}|-1}}
       c_{\vec{\rind},\vec{\sind}}^{\vec{\tind}+\vec{e}_1, (\varphi,\psi)}
       \spair{\vec{\tind}+\vec{e}_1}{\vec{a}\ssha_{(\varphi,\psi)}\vec{b}}
       \quad (\text{by Eq. (\mref{eq:coefzero3})})
  \\
  &=&  \sum_{(\varphi,\psi)\in\indI_{k,\ell}} \sum_{\substack{\vec{\tind}\in \:\ZZ_{\geq 1}^{k+\ell}\\
       |\vec{\tind}|=|\vec{\rind}|+|\vec{\sind}|}}
       c_{\vec{\rind},\vec{\sind}}^{\vec{\tind},(\varphi,\psi)}
       \spair{\vec{\tind}}{\vec{a}\ssha_{(\varphi,\psi)}\vec{b}}.
\end{eqnarray*}
}

\noindent {\bf Case 4. $\rind_1\geq 2$ and $\sind_1=1$.} The proof
for this case is similar to that for Case 3. \qed \medskip

\section{Appendix: a shuffle formulation of the main formula}
\mlabel{sec:appen}
The main body of the paper does not depend on this  Appendix.
 Here we give another formulation of Theorem~\mref{thm:qshsh} in terms of
  shuffles of permutations for those who are interested in a more precise
  connection between the main formula and the shuffle product.

Let integers $k,\ell\geq 1$ be given. Let
\begin{equation}
\begin{aligned}
S(k,\ell):&=\big\{ \sigma\in \Sigma_{k+\ell}\ |\
\sigma^{-1}(1)<\cdots < \sigma^{-1}(k), \sigma^{-1}(k+1)<\cdots
<\sigma^{-1}(k+\ell) \big\}
\\
& = \Big\{\sigma\in \Sigma_{k+\ell}\ \Big |\ \begin{array}{l}
\text{ if } 1\leq \sigma(i)<\sigma(j)\leq k \\
\text{ or } k+1\leq \sigma(i)<\sigma(j) \leq k+\ell, \end{array}
\text{ then } i<j \Big\}.
\end{aligned}
\mlabel{eq:shf}
\end{equation}
be the set of $(k,\ell)$-shuffles.

To state the shuffle form of our main formula we need the following
notations.
Define
$$
\vep_\sigma: [k+\ell]\to \{\pm 1\}, \quad \vep_\sigma(i)=\left\{\begin{array}{ll} 1,& 1\leq \sigma(i)\leq k,\\
    -1, & k+1\leq \sigma(i)\leq k+\ell. \end{array} \right .
$$
Let $\vec{\rind}=(\rind_1,\cdots,\rind_k)\in \ZZ_{\geq 1}^k$ and
$\vec{\sind}=(\sind_1,\cdots,\sind_\ell)\in \ZZ_{\geq 1}^\ell$.
Denote
$$
\vec{\kappa}=(\kappa_1,\cdots,\kappa_{k+\ell}):=(\rind_1,\cdots,\rind_k,\sind_1,\cdots,\sind_\ell).
% \text{ where }
% \kappa_i=\left\{\begin{array}{ll} \rind_i, & 1\leq i\leq k, \\
%    \sind_{i-k}, & \ell+1\leq i\leq k+\ell. \end{array} \right .
$$ Let $\vec{a}\in\sg^k$ and $\vec{b}\in\sg^\ell$. Denote
$$\vec{\gamma}=(a_1,\cdots, a_k,b_1,\cdots, b_\ell). $$ For $\sigma\in
S(k,\ell)$ we denote
$$\vec{a}\ssha_\sigma\vec{b}=(\gamma_{\sigma(1)},\cdots \gamma_{\sigma(k+\ell)}).$$

We have the following equivalent form of Theorem~\mref{thm:qshsh}.

\begin{theorem} %{\bf (Shuffle Form of the Main Theorem)}
Let $\sg$ be a set and let $\calh\bsh(\zsg{\sg})=(\calh(\zsg{\sg}), \qsshab)$ be as defined by Eq.~(\mref{eq:shtransb}).
Then for $\spair{\vec{\rind}}{\vec{a}}\in \zsg{\sg}^k$ and
$\spair{\vec{\sind}}{\vec{b}}\in\zsg{\sg}^\ell$ in $\calh\bsh(\zsg{\sg})$, we have
$$ \spair{\vec{\rind}}{\vec{a}} \qsshab \spair{\vec{\sind}}{\vec{b}} =
\sum_{\substack{
\sigma\in S(k,\ell),\\
\vec{\tind}\in \ZZ_{\geq 1}^{k+\ell}, |\vec{\tind}|=
    |\vec{\rind}|+|\vec{\sind}|}} \Big(\prod_{i=1}^{k+\ell}
    \binc{\tind_i-1}{\kappa_{\sigma(i)}-1 -\frac{1}{2}(1-\vep_\sigma(i)\vep_\sigma(i-1))
    \sum\limits_{j=1}^{i-1}(\tind_j-\kappa_{\sigma(j)})}\Big)
    \spair{\vec{\tind}}{\vec{a}\ssha_{\sigma}\vec{b}}
$$
with the convention that $\vep_\sigma(0)=\vep_\sigma(1)$.
\mlabel{thm:mains}
\end{theorem}
\begin{proof}

Let $\indI_{k,\ell}$ be as defined in Eq.~(\mref{eq:ind}). We have
the bijection between $S(k,\ell)$ and $\indI_{k,\ell}$ given by
\begin{equation}
\sigma^{-1}(j):=\sigma^{-1}_{\varphi,\psi}(j)= \left\{
       \begin{array}{ll} \varphi(j) & \text{ if } 1\leq j\leq k,
                       \\
                       \psi(j-k) & \text{ if } k+1\leq j\leq k+\ell.
       \end{array}
\right.
\mlabel{eq:sigvarpsi}
\end{equation}
That is,
$$
\sigma(i):=\sigma_{\varphi,\psi}(i)= \left\{
       \begin{array}{ll} \varphi^{-1}(i) & \text{ if } i\in \im(\varphi),
                       \\
      k+\psi^{-1}(i) & \text{ if } i\in \im(\psi).
       \end{array}
\right.
$$
Thus we have
\begin{equation}
\kappa_{\sigma(i)}= \left \{ \begin{array}{ll}
\kappa_{\varphi^{-1}(i)}, & i\in \im (\varphi) \\
\kappa_{k+\psi^{-1}(i)}, & i\in \im(\psi)\end{array} \right .
= \left \{\begin{array}{ll} \rind_{\varphi^{-1}(i)}, & i\in \im(\varphi) \\
    \sind_{\psi^{-1}(i)}, & i\in \im(\psi) \end{array} \right .
    = h_{(\varphi,\psi),i}
\mlabel{eq:mumn}
\end{equation}
and
\begin{equation} (\vec{a}\ssha_{\sigma}\vec{b})_i=\gamma_{\sigma(i)}= \left \{ \begin{array}{ll}
\gamma_{\varphi^{-1}(i)}, & i\in \im (\varphi) \\
\gamma_{k+\psi^{-1}(i)}, & i\in \im(\psi)\end{array} \right .
= \left \{\begin{array}{ll} a_{\varphi^{-1}(i)}, & i\in \im(\varphi) \\
    b_{\psi^{-1}(i)}, & i\in \im(\psi) \end{array} \right .
    = (\vec{a}\ssha_{(\varphi,\psi)}\vec{b})_i. \mlabel{eq:mulindeq}
\end{equation} By Eq.~(\mref{eq:mulindeq}) we have
\begin{equation}
\vec{a}\ssha_{\sigma}\vec{b}=\vec{a}\ssha_{(\varphi,\psi)}\vec{b}.
\mlabel{eq:vecab}
\end{equation}

Let $\vep_{\varphi,\psi}$ be the function $[k+\ell]\to \{1,-1\}$
defined in Eq.~(\mref{eq:vep1}). Then for
$\sigma=\sigma_{\varphi,\psi}$,
$$\vep_{\sigma}(i)=1 \Leftrightarrow \sigma(i)\in [k]
\Leftrightarrow i=\sigma^{-1}(j), j\in [k] \Leftrightarrow
i=\varphi(j), j\in [k] \Leftrightarrow i\in \im (\varphi)
\Leftrightarrow \vep_{\varphi,\psi}(i)=1.$$ So we have
\begin{equation}
\vep_{\sigma}(i)=\vep_{\varphi,\psi}(i), \qquad 1\leq i\leq k+\ell.
\mlabel{vepsigma}
\end{equation}
Now our theorem follows from Eq.~(\mref{eq:coef-re-def}),
(\mref{eq:mumn}), (\mref{eq:vecab}), (\mref{vepsigma}) and Theorem
\mref{thm:qshsh}.
\end{proof}

%=======================================================================================
%========================================================================================
%========================================================================================
%\addcontentsline{toc}{section}{\numberline {}References}

\end{document}